\documentclass[thmsa,12pt]{article}
\normalfont\sffamily
\usepackage{amsfonts}
\usepackage{amssymb}
\usepackage[active]{srcltx}
\usepackage{color}
\usepackage[x11names]{xcolor}
\usepackage{amsmath}
\usepackage{amsthm}
\usepackage{bm} % for bold Greek symbols
\usepackage{graphicx}
\usepackage{textcomp}
\usepackage{a4wide}
\usepackage{chicago}
\usepackage{multirow}
\usepackage{diagbox}
\usepackage{color}
\usepackage[format=hang]{caption}
\usepackage{subcaption}
\usepackage{url}
\usepackage{xfrac}

\usepackage[normalem]{ulem}
\usepackage{comment}
\usepackage{pxfonts}
\usepackage[margin=3cm]{geometry} % margin wird weiter unten wieder überschrieben
\usepackage{titlesec}

\usepackage{setspace}
\usepackage{footmisc}
\usepackage{accents}
\usepackage{appendix}

\usepackage{array,booktabs,ragged2e}
\newcolumntype{R}[1]{>{\RaggedRight}p{#1}}

\usepackage{cellspace}
\usepackage{hyperref}
\newtheorem{proposition}{Proposition}

\def \be {\begin{equation}}
\def \ee {\end{equation}}
\newcommand{\proj}[1]{\!\!\downarrow_{#1}}
\newcommand{\lift}[1]{\!\!\uparrow^{#1}}
\DeclareMathOperator*{\argmin}{arg\,min}
\DeclareMathOperator*{\argmax}{arg\,max}

\begin{document}

\def\title #1{\begin{center}
{\Large {\sc #1}}
\end{center}}
\def\author #1{\begin{center} {#1}
\end{center}}

\setstretch{1.1}

\begin{titlepage}
    \phantomsection \label{Titlepage}
    \addcontentsline{toc}{section}{Title page}

\renewcommand{\thefootnote}{\fnsymbol{footnote}}\addtocounter{footnote}{1}
\title{\sc 
Weighted Committee Games \\ \medskip \  }

\author{Sascha Kurz %[Corresponding Author]
	\\ {\small Dept.\ of Mathematics, University of Bayreuth, 95440 Bayreuth, Germany\\email: sascha.kurz@uni-bayreuth.de}}

\author{Alexander Mayer \\ {\small Dept.\ of Economics, University of Bayreuth, 95440 Bayreuth, Germany\\email: alexander.mayer@uni-bayreuth.de}}

\author{Stefan Napel \\ {\small Dept.\ of Economics, University of Bayreuth, 95440 Bayreuth, Germany\\ email: stefan.napel@uni-bayreuth.de}}

\vspace{0.3cm}

\begin{center} {\tt %This draft: 
\hspace{-1.9em} October 31, 2019 %\today
} 
\end{center}

\vspace{0.3cm}
%\vspace{4.01cm}

\begin{center} {\bf {\sc Abstract}} \end{center}
{\small 
Many binary collective choice situations can be described as {weighted simple voting games}. We introduce 
{weighted committee games} to model decisions on an arbitrary number of alternatives in analogous fashion. 
We compare the effect of different voting weights (shareholdings, party seats, etc.) under plurality, Borda, Copeland, and antiplurality rule.
The number and geometry of weight equivalence classes differ widely across the rules. 
Decisions can be much more sensitive to weights in Borda committees than (anti-)plurality or Copeland ones. % rule. 
} 

\vspace{0.2cm}

\begin{description}
{\small
\item[Keywords:] group decisions and negotiations $\cdot$ weighted voting $\cdot$ simple games $\cdot$ scoring rules~$\cdot$ majority rule
%voting games $\cdot$ weighted voting
%$\cdot$ geometry of voting $\cdot$ voting power
%$\cdot$ Borda rule~$\cdot$ Copeland rule  $\cdot$ plurality $\cdot$ antiplurality 
%\item[JEL codes:] C71 $\cdot$ D71 $\cdot$ C63 
% C63: Computational Techniques, Simulation Modeling
% C71: Cooperative Games
% D71: Social Choice, Clubs, Committees, Associations
%\item[Conflict of interest:] The authors declare that they have no conflict of interest.
}
\end{description}

%Highlights:
%- Investors or parties command different vote shares in group choice from m?2 options
%- Many seemingly distinct voting weight distributions imply identical group choices 
%- Weight equivalence classes differ widely for plurality rule, Borda, and Copeland
%- ILP yields minimal representations for all plurality, Borda, and Copeland games 
%- Geometric illustration of which weight change can affect 3-player committee decisions

\vspace{2cm}

\vfill
\noindent {\footnotesize We are grateful to two anonymous referees for their constructive comments. This paper has also benefited from suggestions made by several colleagues and conference audiences, especially H.~Nurmi, N.~Maaser, J.~Sobel, and W.~Zwicker.} 

%\noindent {\footnotesize We are indebted to Hannu Nurmi for drawing our attention to a lacuna 
%that this work addresses. We are also grateful to him, Steven Brams, Dan Felsenthal, Bernard Grofman, Nicola Maaser, Joel Sobel, and William Zwicker for feedback on earlier drafts, as well as to various seminar and conference audiences. }% in Bamberg, Bayreuth, Berlin, Bremen, Dagstuhl, Delmenhorst, Freiburg, Graz, Hagen, Hamburg, Hanover, Leipzig, Moscow, Munich, Rome, Seoul, and Turku.}

\end{titlepage}

\addtocounter{footnote}{-1}
%\newpage
%\pagenumbering{arabic}
%\addtocounter{footnote}{-2}

%%%%%%%%%%%%%%%%%%%%%%%%%%%%%%%%%%%%%%%%%%%%%%%%%%%%%%%%%%%%%%%%%%%%%%%%%
\setstretch{1.26} % war in ET-Version noch 1.2 

\pagenumbering{arabic}

%%%%%%%%%%%%%%%%%%%%
%%%%%%%%%%%%%%%%%%%%
%%%%%%%%%%%%%%%%%%%% INTRODUCTION
%%%%%%%%%%%%%%%%%%%%
%%%%%%%%%%%%%%%%%%%%

\section{Introduction}\label{sec:Introduction}

Consider a corporation of three stockholders owning 6, 5, and 2~million shares. 
As collective decisions between CEO candidates, new business proposals, etc.\ are taken, votes are weighted by the respective shareholdings. 
One may wonder: do the resulting choices differ, ceteris paribus, from those if the three had equal votes? 
Or, say, from outcomes for a $(48\%, 24\%, 28\%)$ distribution? 
The answer affects incentives to participate in a capital increase or to invest in the first place. Similar questions arise in international institutions that weight votes by financial contributions or population sizes (IMF, World Bank, EU~Council) or when parties cast bloc votes in parliament. 

For binary `yes'-or-`no' decisions and a simple majority requirement, the above distributions of weight are equivalent: 
any two shareholders jointly meet the majority threshold of 50\%; the different weights induce the same winning coalitions, hence the same binary voting game. 
This extends to all distributions such that each of three players wields positive but less than half of total weight. A large literature on \emph{simple voting games} has formalized related results.

However, things differ and much less is known if players choose from three or more options. 
For instance, the above shareholders may use plurality rule to decide between three CEO candidates. 
Then investor 1 is decisive whenever 2 and 3 fail to agree: his or her favorite candidate wins with a tally of 6\,:\,5\,:\,2, 
11\,:\,2,  or 8\,:\,5~million. 
Identical plurality winners would result for $(48\%, 24\%, 28\%)$, but ties can arise and yield different decisions for equal votes.
The former weights define the same \emph{committee game}, characterized by $n$~players, $m$~alternatives, and a mapping from $n$-tuples of strict preferences to a winning alternative; equal weights create a different game.

This paper extends several of the~-- unfortunately, already difficult~-- questions that the literature has addressed for binary simple games to committee games. 
Numerical tractability falls sharply in $m$ and $n$ but first answers are feasible.
Our main interest lies in equivalence classes of weight distributions for a given group decision rule. 
One such class, for instance, comprises all weight vectors that induce the same plurality winners as $(6,5,2)$ for decisions on three alternatives.
We investigate three scoring rules (plurality, Borda, antiplurality) and one Condorcet rule (Copeland). 
The respective weight classes can be described by linear inequalities. 
Monotonicity properties \cite{Felsenthal/Nurmi:2017}, degrees of manipulability \cite{Aleskerov/Kurbanov:1999}, or strategic voting equilibria (\citeNP{Myerson/Weber:1993}; %\shortciteNP{Buenrostro/Dhillon/Vida:2013}
\citeNP{Bouton:2013}) for one class member directly apply to all. 

Inspired by the ongoing quest of characterizing and counting all simple voting games that admit a weighted representation,\footnote{See the monograph by \citeN{Taylor/Zwicker:1999} and, e.g., \citeN{Kurz:2012:representation}, \citeN{Houy/Zwicker:2014}, or \shortciteN{Freixas/Freixas/Kurz:2017}. The number of weighted voting games is still unknown for $n>9$.} 
we try to determine the number of distinct weighted committee games involving either plurality, Borda, antiplurality, or Copeland rule for $m\ge 3$ alternatives. 
We propose a decisive test, based on integer linear programming, for whether a given choice rule is representable as a weighted instance of a given scoring rule. 
The superexponential number of  $m^{({m!}\hspace{1pt}^n)}$ distinct mappings from $n$-tuples of preferences to $m$ potential winners makes the analysis computationally difficult. Still, minimal representations are provided for all 51 Borda committee games with $n=m=3$, all plurality and antiplurality committees with $n\le 4$, and all Copeland committees with $n\le 6$ players in the appendix; others are available upon request.
Complete enumeration and lists of games can be useful, e.g., to solve the `inverse problem' of finding a voting game that best achieves a given goal \cite{Kurz:2012} or to obtain sharp bounds on numbers of players and alternatives that permit a certain monotonicity violation, voting paradox, etc. 

The extent to which different voting weights make a real or only a superficial difference has practical relevance. 
For example, financial drawing rights and thereby voting weights among the 24 Directors of the International Monetary Fund's Executive Board were reformed in 2016. Has the vote change been purely cosmetic or is it possible that the agreed weight increases for emerging market economies affect future decisions, such as the choice of the next IMF Managing Director? 
The selection process for the latter has been reformed too and involves a shortlist of three candidates that is compiled ``taking into account the Fund's weighted voting system''.\footnote{Cf.~{www.imf.org/en/About/Factsheets/Managing-Director-Selection-Process/qandas} (last accessed: August 12, 2019).} 
The winner is chosen by consensus or else determined ``by a majority of the votes cast'' according to the new voting powers. 
The IMF's weighted voting system would be compatible with pairwise majority comparisons of candidates (Copeland rule) as well as plurality votes in the shortlisting and runoff stages; in the former, use of a scoring method such as Borda rule is conceivable too. 
Whether different procedural choices at either stage might have an effect and whether the 2016 vote changes could matter for outcomes both boil down to equivalence or not of weighted committee games. 

We depict the geometry of weighted committees for $n=3$. The illustrations convey a sense of the (non-)robustness of collective decisions with respect to variations in voting weights that reflect, e.g., changes in corporate voting rights or party switching in parliament. Borda rule can be seen to be highly sensitive to weight variations while collective choice by (anti-)plurality or Copeland rule generates only muted incentives for increasing votes to gain influence on decisions.

%%%%%%%%%%%%%%%%%%%%
%%%%%%%%%%%%%%%%%%%%
%%%%%%%%%%%%%%%%%%%% SECTION 2: Preliminaries 
%%%%%%%%%%%%%%%%%%%%
%%%%%%%%%%%%%%%%%%%%

\section{Notation and definitions}\label{sec:notation_and_definitions}
\subsection{Preliminaries}

We consider a finite set $N$ of $n\ge 1$ players such that each $i\in N$ has a strict preference relation $P_i$ over a set $A=\{a_1,\dots,a_m\}$ of $m\ge 2$ alternatives. 
$\mathcal{P}(A)$ denotes the set of all $m!$ strict preference orderings on $A$. 
A \emph{collective choice rule} $\rho\colon  \mathcal{P}(A)^n \to A$ maps each profile $\mathbf{P}=(P_1,\dots,P_n)$ to a winning alternative $a^*=\rho(\mathbf{P})$.  
Though $\rho$ is defined on complete preference profiles $\mathbf{P}$, it might draw on partial preference information only and require players, e.g., to submit just the top elements of  $P_1,\ldots, P_n$ in a plurality vote. Note also that $\rho$ does not specify information on how non-winning alternatives $a\neq a^*$ are ranked relative to each other: we investigate resolute (singleton-valued) choice rules rather than Arrovian social welfare functions.

Rules $\rho$ that treat all voters $i\in N$ symmetrically play a special role in our analysis: suppose profile $\mathbf{P'}=(P_{\pi(1)},  \ldots, P_{\pi(n)})$ results from applying a permutation $\pi\colon N\to N$ to $\mathbf{P}$. 
Then $\rho$ is  \emph{anonymous} if $\rho(\mathbf{P})=\rho(\mathbf{P'})$ for all such $\mathbf{P}$, $\mathbf{P'}$. 
We  write $r$ instead of $\rho$ if we want to highlight that a considered rule is anonymous. 

We focus on four standard voting rules. %with lexicographic tie breaking. 
Their definitions are summarized in Table~\ref{table:rules} where $b_i(a,\mathbf{P}):= |\{a'\in A \ |\ a P_i a'\}|$. 
Antiplurality rule $r^A$, Borda rule $r^B$, and plurality rule $r^P$ are \emph{scoring rules}:
winners can be characterized as maximizers of scores derived from alternatives' positions in $\mathbf{P}$ and a fixed {scoring vector} $\mathbf{s}\in \mathbb{Z}^m$ with $s_1\ge s_2\ge \ldots \ge s_m$. 
Namely, let the fact that alternative $a$ is ranked at the $j$-th highest position in ordering $P_i$ contribute $s_j$ points for $a$, and refer to the sum of all points received as $a$'s \emph{score}. Then score maximization for $\mathbf{s}^B=(m-1, m-2, \ldots,1, 0)$ yields the Borda winner, $\mathbf{s}^P=(1, 0, \ldots, 0, 0)$ the plurality winner, and $\mathbf{s}^A=(0, 0, \ldots, 0, -1)$ or $(1, 1, \ldots, 1, 0)$ the antiplurality winner.
By contrast, {Copeland rule} $r^C$ considers the pairwise majority relation 
$
a\succ_M^\mathbf{P} a' \ :\Leftrightarrow \  \big|\{i\in N\ | \ a P_i a' \} \big| >  \big|\{i\in N\ | \ a' P_i a \}\big|
$ and selects
the alternative that beats the most others according to $\succ_M^\mathbf{P}$.  Clearly, if $\succ_M^\mathbf{P}$ has a unique top~-- the \emph{Condorcet winner}~-- then it will be chosen: that is, 
$r^C$ is a \emph{Condorcet method}.

\renewcommand{\arraystretch}{1.4}
\begin{table}
	\begin{center}
		\begin{tabular}{|l|l|}
			\hline
			Rule & Winning alternative at preference profile $\mathbf{P}$ \\
			\hline 
			Antiplurality %$r^A$ 
			&  $r^A(\mathbf{P})\in \argmin_{a\in A} \big|\{i\in N \ |\ \forall a'\neq a\in A\colon a' P_i a\}\big|$
			\\
			Borda %$r^B$ 
			&  $r^B(\mathbf{P})\in \argmax_{a\in A} \sum_{i\in N} b_i(a,\mathbf{P})$\\ 
			Copeland %$r^C$ 
			& $r^C(\mathbf{P})\in \argmax_{a\in A} \big|\{a'\in A \ |\ a \succ_M^\mathbf{P} a'\}\big|$ \\
			Plurality %$r^P$ 
			&  $r^P(\mathbf{P})\in \argmax_{a\in A} \big|\{i\in N \ |\ \forall a'\neq a\in A\colon a P_i a'\}\big|$ \\
			\hline 
		\end{tabular} 
		\caption{Investigated voting rules \label{table:rules}}
	\end{center}
\end{table}

We impose {lexicographic tie breaking}. 
This has computational advantages over working with set-valued choices and entails no loss provided we consider all $\mathbf{P}\in \mathcal{P}(A)^n$: 
the set of alternatives tied at $\mathbf{P}$ is fully determined by $a^*=r(\mathbf{P})$ and the respective winners $a^{**}, a^{***}, \ldots$ at profiles $\mathbf{P'},\mathbf{P''}, \ldots$ that swap $a^*$ with alternatives $a', a'', \ldots$ %in $\mathbf{P}$ 
that might be tied with $a^*$ at $\mathbf{P}$.\footnote{Given $r(\mathbf{P})=b$, for example, a tie with $a$ can directly be ruled out; 
	one sees if $b$ was tied with $c$ by checking whether $r(\mathbf{P'})=c$ or~$b$ where $\mathbf{P'}$ only swaps $b$'s and $c$'s position in every player's ranking $P_i$.} 
Our resolute rules $r^A$, $r^B$, $r^C$, $r^P$ and their set-valued versions are hence in one-to-one correspondence and exhibit identical structural equivalences. The same applies to uniform random tie breaking.

We call the combination $(N, A, \rho)$ of a set of voters, a set of alternatives and a collective choice rule a \emph{committee game} or just a \emph{committee}. 
Several special cases have previously received attention in the literature. 
 
Most prominently, committees $(N, A, \rho)$ with binary $A=\{0,1\}$ and surjective, monotonic $\rho$ are known as \emph{simple voting games}. Following \citeN[Ch.~10]{vonNeumann/Morgenstern:1953}, these games are usually described as a pair $(N,v)$ with %where $\rho$ is replaced by {function} 
$v\colon 2^N \to \{0,1\}$ and $v(S)=1$ when $1\,P_i \,0$ for all $i\in S$ implies $\rho(\mathbf{P})=1$. Sets $S\subseteq N$ with $v(S)=1$ are known as \emph{winning coalitions}. 
$(N,v)$ is \emph{weighted} and called \emph{weighted voting game} if there exists a non-negative %, non-degenerate 
vector $\mathbf{w}=(w_1, \ldots, w_n)$ of weights and a positive quota $q$ such that $v(S)=1$ if and only if $ %\Leftrightarrow 
\sum_{i\in S} w_i \ge q$.\footnote{
Existence is guaranteed only for $n\le 3$. Games that are not weighted arise, e.g., in the Legislative Council of Hong Kong \cite{Cheung/Ng:2014} and the EU~Council \cite{Kurz/Napel:2016}.} 
Pair $(q; \mathbf{w})$ is a \emph{(weighted) representation} of $(N,v)$ and one writes % denotes the respective game by $[q; \mathbf{w}]$, i.e.,
$(N,v)=[q; \mathbf{w}]$. 
It is no restriction to focus on integers: given $q\in \mathbb{R}_{++}$, $\mathbf{w}\in \mathbb{R}^n_+$ there always exist $q'\in \mathbb{N}$, $\mathbf{w'}\in \mathbb{N}^n_0$ such that $[q; \mathbf{w}]=[q'; \mathbf{w'}]$. 

Other special cases include voting with multiple levels of approval, such as \emph{ternary voting games} (\citeNP{Felsenthal/Machover:1997}), %\shortciteN{Tchantcho/DiffoLambo/Pongou/Engoulou:2008} and \citeN{Parker:2012} have considered \emph{ternary voting games} with the %individual 
%option to support a proposal, to abstain, or to reject it. 
\emph{quaternary voting games} (\citeNP{Laruelle/Valenciano:2012}), and \emph{$(j,k)$-games} %in which each player goes for one of $j$ linearly ordered options and every resulting $j$-partition of player set $N$ is mapped to one of $k\le j$ ordered outcomes 
(\citeNP{Hsiao/Raghavan:1993}; \citeNP{Freixas/Zwicker:2003}).  %, \citeyearNP{Freixas/Zwicker:2009}). 
Plurality committees (defined below) have featured in the multicandidate voting frameworks of %\citeN{Bolger:1986}
\citeN{Bolger:1983}, \shortciteN{Amer/Carreras/Magana:1998b} and Monroy and Fern\'{a}ndez~\citeyear{Monroy/Fernandez:2009,Monroy/Fernandez:2011}
as \emph{simple plurality games} and \emph{relative majority $r$-games}.\footnote{ 
Their analysis focused on power indices. 
We here show that there are only 36 distinct simple plurality games with four players and so at most 36 different distributions of power can arise.}

\subsection{Weighted committee games}\label{subsec:weighted_committee_games}

Many committee games that model real collective decision making involve a {non-anonymous} rule $\rho$. 
For instance, an anonymous decision rule $r$ may apply at the level of shareholdings, IMF~drawing rights, etc.\ rather than that of individual voters; or players $i\in N$ are well-disciplined parties with different numbers of seats.\footnote{
If designated members enjoy procedural privileges, this may~-- but need not~-- be equivalent to asymmetric weights. Veto power of permanent members in the UN~Security Council, e.g., translates into $[39; 7,7,7,7,7, 1, \ldots, 1]$ for $m=2$.  } 
The corresponding rule $\rho$ can be viewed as the \emph{combination of an anonymous collective choice rule $r$ with integer voting weights $w_1, \ldots, w_n$} attached to the players.

With $r$ denoting the entire family of mappings from $n$-tuples of linear orders over $A$ to winners $a^*\in A$ under the considered rule, 
we define $r|\mathbf{w}\colon \mathcal{P}(A)^{w_\Sigma}\to A$ by
\be
r|\mathbf{w}(\mathbf{P}):=r(\underbrace{P_1, \ldots, P_1}_{\textstyle w_1 \text{ times}},\  \underbrace{P_2, \ldots, P_2}_{\textstyle w_2 \text{ times}},\ \ldots,\  \underbrace{P_n, \ldots, P_n}_{\textstyle w_n \text{ times}}) \label{eq:rw_defined}
\ee
for a given anonymous rule $r$ and a non-negative, non-degenerate weight vector $\mathbf{w}=(w_1, \ldots, w_n)\in \mathbb{N}_0^n$ with $w_\Sigma:=\sum_{i=1}^n w_i>0$. 
In the degenerate case $\mathbf{w}=(0, \ldots, 0)$, let $r|\mathbf{0}(\mathbf{P})\equiv a_1$.

We say a committee game $(N, A, \rho)$ is \emph{$r$-weighted} for a given rule $r$ if there exists a weight vector $\mathbf{w}=(w_1, \ldots, w_n)\in \mathbb{N}_0^n$ such that $(N, A, \rho)=(N,A,r|\mathbf{w})$, i.e., 
\be
\rho(\mathbf{P})=r|\mathbf{w}(\mathbf{P}) \text{ for all }\mathbf{P}=(P_1, \ldots, P_n)\in \mathcal{P}(A)^n.
\ee
$(N,A,r,\mathbf{w})$ is a \emph{(weighted) representation} of $(N, A, \rho)$ and
we also denote this game by $[N,A,r,\mathbf{w}]$. %, i.e., $(N,A,r|\mathbf{w})=[N,A,r,\mathbf{w}]$.
$[N,A,r^A,\mathbf{w}]$, $[N,A,r^B,\mathbf{w}]$, $[N,A,r^C,\mathbf{w}]$ and $[N,A,r^P,\mathbf{w}]$ are referred to as \emph{antiplurality}, \emph{Borda}, \emph{Copeland} and \emph{plurality committees} respectively.
\renewcommand{\arraystretch}{1.2}
\begin{table}
\begin{center}
\begin{tabular}{|c|c|c|c|}
\hline
$P_1$ & $P_2$ & $P_3$ & $P_4$ \\
 \hline 
$ d $ & $ b $ & $ c $ & $ c $ \\
$ e $ & $ c $ & $ e $ & $ b $ \\
$ b $ & $ e $ & $ a $ & $ a $ \\
$ a $ & $ a $ & $ d $ & $ d $ \\
$ c $ & $ d $ & $ b $ & $ e $ \\
\hline 
\end{tabular} 
\quad $\Rightarrow$ \quad
\begin{tabular}{cl}
$r^A|\mathbf{w}(\mathbf{P})=a$ & ($a$ has min.\ %bottom ranks 
negative votes 0) \\
$r^B|\mathbf{w}(\mathbf{P})=b$ & ($b$ has max.\ Borda score 28)\\
$r^C|\mathbf{w}(\mathbf{P})=c$ & ($c$ has max.\ pairwise wins 3) \\ 
$r^P|\mathbf{w}(\mathbf{P})=d$ & ($d$ has max.\ plurality tally 5) 
\end{tabular} 
\caption{Choices for preference profile $\mathbf{P}$ when $\mathbf{w}=(5,3,2,2)$\label{table:example}}
\end{center}
\end{table}
Such committees typically differ for $m>2$, as illustrated in Table~\ref{table:example}: 
the winner from $A=\{a,b,c,d,e\}$ at $\mathbf{P}$ all depends on the voting rule in use. Neither this observation nor below structural findings depend on whether $\mathbf{P}$ reflects sincere or strategic preference statements.

By definition (\ref{eq:rw_defined}) uniform weights $\mathbf{w}=(1,\ldots, 1)$ reduce any collective choice rule $r|\mathbf{w}$ to $r$. 
It follows that monotonicity, consistency, and other properties that are satisfied (violated) by a given rule $r$ are also satisfied (violated) by $r|\mathbf{w}$. 
For instance, the axiomatic characterizations by \citeN{Young:1975} of set-valued versions of general scoring rules, \citeN{Henriet:1985} of a 
Copeland variation, or \citeN{Kurihara:2018} of antiplurality continue to apply. 
Non-uniform weights $\mathbf{w}$ essentially impose the domain restriction that $w_1$, $w_2$, $\ldots$ individuals have identical preferences.
At a given profile $\mathbf{P}$, a group~$i$ of voters can be worse off under $r|\mathbf{w'}$ than $r|\mathbf{w}$ for $w'_i>w_i$ and $w'_j=w_j$, $j\neq i$, if and only if $r$ suffers from the so-called \emph{no show paradox}. Condorcet methods like Copeland rule are known to do so (see \citeNP{Moulin:1988}), while scoring rules do not.
Our question here is related but a broader one: for given $r$, which weight variations $\mathbf{w'}\neq \mathbf{w}$ can or cannot induce choice changes compared to $r|\mathbf{w}$?

\section{Equivalence classes of weighted committee games}\label{sec:equivalence_of_committees}

\subsection{Equivalence of committee games}

Two $r$-weighted committee games $(N, A,r|\mathbf{w})$ and $(N', A',r|\mathbf{w'})$ are \emph{structurally equivalent} or \emph{equivalent up to isomorphism} if 
\begin{equation}
\big\{a_j\,P_i\, a_k \Leftrightarrow \tilde\pi(a_j)\, P_{\pi(i)}'\, \tilde\pi(a_k)  \big\}\, \Rightarrow\  \tilde\pi\big(r|\mathbf{w}(\mathbf{P})\big)=r|\mathbf{w'}(\mathbf{P'})
\end{equation}
for bijections $\pi\colon N\to N'$ and $\tilde\pi\colon A\to A'$ that map every profile $\mathbf{P}$ of preferences $P_i$ over $A$ to a relabeled profile $\mathbf{P'}$ of preferences $P'_{\pi(i)}$ over $A'$.\footnote{Analogous equivalences apply to $(N, A,r|\mathbf{w})$ and $(N', A',r'|\mathbf{w'})$ when $r\neq r'$ or general committees $(N, A,\rho)$ and $(N', A',\rho')$.  
The respective considerations could, in principle, also be extended to Arrovian \emph{social welfare functions}, which map each profile %$\mathbf{P}$ 
of individual preferences to a collective preference $P^*\in \mathcal{P}(A)$ rather than a winning alternative $a^*\in A$. We leave detailed explorations to future research.} 
For instance, Copeland committees with $N=N'=\{1,2,3\}$, $A=A'$, $\mathbf{w}=(3,1,1)$ and $\mathbf{w'}=(1,3,1)$
have quite different attractiveness to player~1 but the decision environment is structurally the same: there is a dictator player whose most-preferred alternative wins and two null players whose preferences do not affect the outcome. 

A given distribution $\mathbf{w}\in \mathbb{N}_0^n$ fixes $n$ and we can write $(r,\mathbf{w})\sim_{m} (r,\mathbf{w'})$ if $r$-committee games 
with $m$~alternatives are structurally equivalent for weights $\mathbf{w}$ and $\mathbf{w'}$.
%The %equivalence 
Relation $\sim_{m}$ and some 
$\mathbf{\bar w} \in \mathbb{N}_0^n$ with $\bar w_1 \ge \bar w_2\ge \ldots \ge \bar w_n$ %serving as index 
define the \emph{equivalence class}
\be
\mathcal{E}^{r}_{\mathbf{\bar w},m}:=\big\{\mathbf{w} \in \mathbb{N}_0^n \ | \ (r,\mathbf{w})\sim_{m}  (r,\mathbf{\bar w})\big\}.
\ee
If rule~$r$ is in use for deciding between $m$ alternatives, then all weight distributions $\mathbf{w}, \mathbf{w'}\in \mathcal{E}^{r}_{\mathbf{\bar w},m}$ come with identical monotonicity properties, voting paradoxes, manipulation incentives, implementation possibilities, strategic equilibria,
etc.

\subsection{Illustration}\label{subsec:illustration}

As an example, consider Borda rule for $m=3$ and reference weights $\mathbf{\bar w}=(5,2,1)$ that reflect a given seat distribution in a council, voting stocks held by shareholders, etc. 
We focus on the subset $\mathcal{E}^{B}_{(5,2,1),3}\subset \mathcal{E}^{r^B}_{(5,2,1),3}$ 
of alternative distributions $\mathbf{w}$ with $w_1\ge w_2\ge w_3$. 
Two linear inequalities are implied by $r^B|(5,2,1)=r^B|\mathbf{w}$ for each profile $\mathbf{P}\in\mathcal{P}(A)^3$.
For instance, writing $abc$ in abbreviation of $a P_i b P_i c$, profile $\mathbf{P}=(cab,bac,abc)$
gives rise to a total Borda score $\bar w_1\cdot 1 + \bar w_2\cdot 1 + \bar w_3 \cdot 2=9$ for alternative $a$: it is ranked top by player~3, middle by players~1 and 2. The corresponding scores for $b$ and $c$ are $2\bar w_2+\bar w_3=5$ and $2\bar w_1=10$. Hence $r^B|\mathbf{\bar w}(\mathbf{P})=c$ and any allocation  $\mathbf{w}$ of seats, shares, etc.\ that is equivalent to $\mathbf{\bar w}$ must ensure that the respective Borda score $2w_1$ of (lexicograpically maximal) $c$ strictly exceeds $a$'s and $b$'s scores:
\begin{align}
&\text{(I) }\ 2w_1 > w_1+w_2+2w_3  & \text{and} & & &\text{(II) }\
2w_1 > 2w_2+w_3.\notag
\intertext{$\mathbf{P'}=(cab,abc,bac)$ makes $a$ the winner. 
Its score must not be smaller than $b$'s and $c$'s:}
&\text{(III) }\ w_1+2w_2+w_3 \ge w_2+2w_3  & \text{and} & & &\text{(IV) }\
w_1+2w_2+w_3\ge 2w_1.\notag
\intertext{Wins by $a$ and $b$ for $\mathbf{P''}=(abc,bca,bac)$ and $\mathbf{P'''}=(abc,bca,bca)$ 
similarly imply:}
&\text{(V) }\ 2w_1+w_3 \ge w_1+2w_2+2w_3  & \text{and} & & &\text{(VI) }\
2w_1+w_3 \ge w_2\notag \\
&\text{(VII) }\ w_1+ 2w_2+2w_3 > 2w_1  & \text{and} & & &\text{(VIII) }\
w_1+ 2w_2+2w_3 \ge w_2+w_3.\notag
\end{align}
Condition~(VIII) is trivially satisfied. (IV) and (V) imply $w_1 = 2w_2+w_3$. This makes (I) equivalent to $w_2>w_3$ and (VII) to $w_3>0$. Combining $w_1 = 2w_2+w_3$ and $w_2>w_3>0$ also verifies (II), (III) and (VI).
The 212 remaining profiles $\mathbf{P}\in \mathcal{P}(A)^3$ turn out not to impose additional constraints. Hence 
\be
\mathbf{w}\in \mathcal{E}^{B}_{(5,2,1),3}=\big\{(2w_2+w_3,w_2,w_3)\in\mathbb{N}_0^3\,:\,w_2> w_3>0 
\big\}
\ee
contains all weight distributions $w_1\ge w_2\ge w_3$ that imply Borda choices identical to $\mathbf{\bar w}=(5,2,1)$ for all preference profiles over three options.
The full class $\mathcal{E}^{r^B}_{(5,2,1),3}$ follows by permuting the distributions in $\mathcal{E}^{B}_{(5,2,1),3}$. 
Other classes, such as $\mathcal{E}^{r^B}_{(1,1,1),3}$, $\mathcal{E}^{r^B}_{(2,1,1),3}$, etc., are characterized by analogous inequalities. 

\subsection{Relation between equivalence classes}

\label{subsec_theoretical}

The number of distinct mappings from preference profiles to outcomes is large but finite for given $n$ and $m$.
$\mathbb{N}_0^n$ hence is partitioned into a finite collection
$
\big\{\mathcal{E}^{r}_{\mathbf{\bar w_1},m}, \mathcal{E}^{r}_{\mathbf{\bar w_2},m}, $ $\ldots, \mathcal{E}^{r}_{\mathbf{\bar w}_{\bm{\xi}},m}\big\}
$
for any given rule $r$. 
We will investigate numerically how the number $\xi$ of elements varies across rules but let us first state some analytical observations.
The first two are obvious: 

\begin{proposition} \label{prop:all_the_same}
The partitions $\big\{\mathcal{E}^{r}_{\mathbf{\bar w_1},2}, \ldots, \mathcal{E}^{r}_{\mathbf{\bar w}_{\bm{\xi}},2}\big\}$ of $\mathbb{N}_0^n$ coincide for $r\in \{r^A, r^B, r^C, r^P\}$.
\end{proposition}

\begin{proposition} \label{prop:wvg_plurality}
Let $A=\{a_1,a_2\}$ and $r\in \{r^A, r^B, r^C, r^P\}$. Then 
$
r|\mathbf{w}(\mathbf{P})=a_1 \Leftrightarrow v(S)=1
$
where $(N,v)=[q;\mathbf{w}]$ with $q=\frac{1}{2}\sum_{i\in N} w_i$ and $S=\{i\in N\ | \ a_1\, P_i\, a_2\}.$
\end{proposition}

\noindent  
It follows that the respective partitions $\big\{\mathcal{E}^{r}_{\mathbf{\bar w_1},2}, 
\ldots, \mathcal{E}^{r}_{\mathbf{\bar w}_{\bm{\xi}},2}\big\}$ of $\mathbb{N}_0^n$ 
coincide with those for weighted voting games with a simple majority quota. Their study and enumeration for $n\le 5$ dates back to \citeN[Ch.~10]{vonNeumann/Morgenstern:1953}. %Also see \citeN{Brams/Fishburn:1996}.

The next observations vary $m$ for fixed $r$; proofs are given in Appendix~A:
\begin{proposition}\label{prop:antiplurality}
	The antiplurality partitions $
	\big\{\mathcal{E}^{r^A}_{\mathbf{\bar w_1},m}, \mathcal{E}^{r^A}_{\mathbf{\bar w_2},m}, 
	\ldots, \mathcal{E}^{r^A}_{\mathbf{\bar w}_{\bm{\xi}},m}\big\}
	$  of $\mathbb{N}_0^n\smallsetminus \{\mathbf{0}\}$ consist of $\xi=n$ equivalence classes identified by weight vectors 
	$\mathbf{\bar w_1}=(1,0,\ldots, 0), 
	\mathbf{\bar w_2}=(1,1,\ldots, 0), 
	\ldots, 
	\mathbf{\bar w_n}=(1,1,\ldots, 1) $
	for all $m\ge n+1$.
\end{proposition}

\begin{proposition} \label{prop:Borda}
	For Borda rule $r^B$ and $m\ge 3$, every weight vector $\mathbf{\tilde w_j}=(j,1,0, \ldots, 0)$ with $j\in \{1, \ldots, m-1\}$ 
	identifies a different class $\mathcal{E}^{r^B}_{\mathbf{\tilde w_j},m}$.
\end{proposition}

\noindent This implies that~-- differently from antiplurality, Copeland, and plurality~-- the number $\xi$ of structurally distinct Borda committee games for given $n\ge 3$ grows in $m$ without bound.

\begin{proposition} \label{prop:Copeland}
The Copeland %rule %$r^C$, the 
partitions
$
\big\{\mathcal{E}^{r^C}_{\mathbf{\bar w_1},m}, \ldots, \mathcal{E}^{r^C}_{\mathbf{\bar w}_{\bm{\xi}},m}\big\}
$ of $\mathbb{N}_0^n$ coincide for all $m\ge 2$.
\end{proposition}

\noindent So $r^C$ extends the known equivalences for binary simple voting games to arbitrarily many options. 
This might feel unsurprising because winners in Copeland committees are selected by binary comparisons. 

However, the conjecture that Prop.~\ref{prop:Copeland} applies to just any Condorcet method is wrong.
Copeland rule is special. 
For instance, \emph{Black rule} selects the Condorcet winner if one exists and otherwise breaks cyclical majorities by Borda scores. 
Weight distributions of $(6,4,3)$ and $(4,4,2)$ are equivalent for $m=2$ and give rise to a cycle $a \succ_M^\mathbf{P} b \succ_M^\mathbf{P} c \succ_M^\mathbf{P} a$ for $\mathbf{P}=(cab,abc,bca)$. 
The Black winner then is $c$ for the former but $a$ for the latter weights; so they are non-equivalent for $m=3$. The same applies to \emph{Kemeny--Young} or \emph{maximum likelihood rule}, which picks the top element of the collective preference ranking $P^*$ that minimizes total pairwise disagreements %(Kemeny distances) 
with all individual rankings in $\mathbf{P}$; or \emph{maximin rule}, where the winner maximizes the minimum support across all pairwise comparisons.

\begin{proposition}\label{prop:plurality}
The %For plurality rule $r^P$, the 
plurality partitions
$
\big\{\mathcal{E}^{r^P}_{\mathbf{\bar w_1},m}, %\mathcal{E}^{r^P}_{\mathbf{\bar w_2},m}, 
\ldots, \mathcal{E}^{r^P}_{\mathbf{\bar w}_{\bm{\xi}},m}\big\}
$ of $\mathbb{N}_0^n$ coincide for all $m\ge n$.

\end{proposition}

%There are at least $m-1$ Borda committees involving $n\ge 2$ players and $m\ge 3$ alternatives: 

%\begin{proposition} \label{prop:Borda}
%For Borda rule $r^B$ and given $m\ge 3$, every weight vector $\mathbf{\tilde w_j}=(j,1,0, \ldots, 0)$ with $j\in \{1, \ldots, m-1\}$ 
%identifies a different %equivalence 
%class $\mathcal{E}^{r^B}_{\mathbf{\tilde w_j},m}$.
%\end{proposition}

%\begin{proof} 
%\noindent Let $k>j$ for otherwise arbitrary $j, k\in \{1, \ldots, m\}$. Consider $A=\{a_1, \ldots, a_m\}$ and any profile $\mathbf{P}\in\mathcal{P}(A)^n$ such that player~1 prefers $a_2$ most and ranks all remaining alternatives lexicographically
%while player~2 ranks $a_2$ in $k$-th position and otherwise agrees with player~1, i.e., suppose 
%$a_2\,P_1\,a_1\,P_1\,a_3\,P_1\,a_4\,%P_1\,
%\ldots\, %P_1\, 
%a_m$ and 
%$a_1\,P_2\,a_3\,P_2\,a_4 %\,P_2
%\,\ldots$ $%\, P_2\, 
%a_k\,P_2\, a_2 \,P_2\, a_{k+1}\,P_2\, a_{k+2}\, %P_2\, 
%\ldots \, %P_2\, 
%a_{m}.$

\section{Identification of weighted committees} \label{sec:identification}

\subsection{Minimal representations and test for weightedness}\label{subsec:minimal_representations}

Given $(N, A, \rho)=(N,A,r|\mathbf{w})$,
we say that $(N, A,r,\mathbf{w})$ has \emph{minimum integer sum} or is a \emph{minimal representation} of $(N,A,\rho)$ if $\sum_{i\in N}w'_i\ge \sum_{i\in N}w_i$ for all representations $(N,A,r,\mathbf{w'})$ of $(N,A,\rho)$.
Games in a given equivalence class $\mathcal{E}^{r}_{\mathbf{\bar w},m}$
usually have a unique minimal representation;\footnote{For $m=2$, minimal representations are unique up to $n=7$ players \cite{Kurz:2012:representation}. Multiplicities for games with larger values of $m$ or $n$ arise but are rare.} 
corresponding weights are a focal choice for $\mathbf{\bar w}$. 
Minimal representations can be more informative or convenient to work with than given weights in applications (cf.\ \citeNP{Freixas/Kaniovski:2014}, \citeNP{Kurz/Napel:2016}). 

Finding minimal representations of all Copeland committees simplifies to finding weighted representations of simple voting games (Prop.~\ref{prop:wvg_plurality} and \ref{prop:Copeland}). Linear programming techniques have proven useful there and can be adapted to scoring rules $r^A$, $r^B$, or $r^P$. 
Namely, consider a scoring rule $r$ based on an arbitrary but fixed scoring vector~$\mathbf{s}$. Write $S_{k}(P_i)\in\mathbb{Z}$ for the unweighted $(s_1, s_2, \ldots, s_m)$-score of alternative $a_k$ derived from its position in ordering $P_i$; for instance, for $m=3$ and $a_3=c$, we have  $S_{3}(P_i)=s_2$ if either $a P_i c P_i b$ or $b P_i c P_i a$. 
Now suppose that $r$ can be combined with integer voting weights so as to induce choice rule~$\rho$. Then~-- denoting the index of the winning alternative at profile $\mathbf{P}$ by $\omega_{\rho}(\mathbf{P})\in\{1,\dots,m\}$, i.e., 
$\rho(\mathbf{P})=a_{\omega_{\rho}(\mathbf{P})}\in A$~-- any solution to the following \emph{integer linear program} yields a minimal representation $(N, A, r, \mathbf{w})$ of $(N, A, \rho)$:
\begin{align}
&  \min_{\mathbf{w}\in\mathbb{N}_0^n} \ \sum_{i=1}^n w_i %\quad\text{s.t.}
\label{eq:ilp_min_sum}\tag{ILP} & \\ \text{s.t.}\quad
&  \sum_{i=1}^n S_{k}(P_i)\cdot w_i \le \sum_{i=1}^n S_{\omega_{\rho}(\mathbf{P})}(P_i)\cdot w_i\,-\,1\quad \forall\mathbf{P}\in\mathcal{P}(A)^n\, \forall 1\le k\le \omega_{\rho}(\mathbf{P})-1,\nonumber\\  
&  \sum_{i=1}^n S_{k}(P_i)\cdot w_i \le \sum_{i=1}^n S_{\omega_{\rho}(\mathbf{P})}(P_i)\cdot w_i\quad \forall\mathbf{P}\in\mathcal{P}(A)^n\, \forall \omega_{\rho}(\mathbf{P})+1\le k\le m.\nonumber
\end{align}

\noindent The case distinction between %scores of non-winning 
alternatives $a_k$ with index $k< \omega_{\rho}(\mathbf{P})$ vs.\ $k> \omega_{\rho}(\mathbf{P})$ reflects the tie breaking assumption. If some (non-minimal) representation $(N, A, r, \mathbf{w'})$ of $(N, A, \rho)$ is %already 
known and $w'_1\ge w'_2\ge \ldots \ge w'_n$ then adding constraints 
$
w_i\ge w_{i+1}, \forall 1\le i\le n-1,
$ 
to (\ref{eq:ilp_min_sum}) accelerates computations. 

If it is not known whether $\rho$ is $r$-weighted, \eqref{eq:ilp_min_sum} provides a decisive \emph{test for $r$-weightedness} for any given scoring rule~$r$ (and $r^C$ by Prop.~\ref{prop:all_the_same} and \ref{prop:Copeland}): the constraints characterize a non-empty compact set if and only if $\rho$ is $r$-weighted. Checking non-emptiness of \eqref{eq:ilp_min_sum}'s constraint set can be done with software (e.g., Gurobi or CPLEX) that determines a weight sum minimizer at little extra effort.

\subsection{Algorithm for identifying all $r$-committees} 
\label{subsec_algorithmic_methods}

In principle, one could find and characterize {all} distinct $r$-committee games for fixed $n$ and $m$ as follows:
loop over all $m^{(m!^n)}$ mappings $\rho\colon \mathcal{P}(A)^n\to A$; conduct above test; 
in case of success, determine a representation $(N,A,r,\mathbf{\bar w})$
and characterize  $\mathcal{E}_{\mathbf{\bar w}, m}^r$ as in Section~\ref{subsec:illustration}; continue until all choice rules $\rho$ have been covered. 

The explosive growth of $m^{(m!^n)}$ prevents a direct implementation of this idea.\footnote{$3^{(3!^3)}=3^{216}> 10^{103}$ already exceeds the estimated number of % $%\approx\! %	10^{80}$ %10^78--$10^{84}$ 
	atoms in the universe.}
However, many mappings can be dropped from consideration. 
If $\rho(\mathbf{P})=a_1$ for one of the $(m-1)!^n$ profiles $\mathbf{P}$ where $a_1$ is unanimously ranked last, for instance, then  $\rho$ cannot be $r$-weighted for $r\in \{r^A, r^B, r^C, r^P\}$. This rules out $m^{(m!^n-1)}$ candidate mappings in one go. Similarly, if weights $\mathbf{w}$ such that $r|\mathbf{w}(\mathbf{P})=a_1$ turn out to be incompatible with $r|\mathbf{w}(\mathbf{P'})=a_2$ for two suitable profiles $\mathbf{P}, \mathbf{P'}$, then all $m^{(m!^n-2)}$ mappings $\rho$ with $\rho(\mathbf{P})=a_1$ and $\rho(\mathbf{P'})=a_2$ can be disregarded at once.
The branch-and-cut algorithm described in Table~\ref{tab:BC_algorithm} operationalizes these considerations.

\renewcommand{\arraystretch}{1.4}

\begin{table}
\centering 
\begin{tabular}{l p{12.5cm}}
\hline 
\multicolumn{2}{l}{\textbf{Branch-and-Cut Algorithm}\vspace{-0.2cm}} 
\\
\multicolumn{2}{l}{Given $n$, $m$ and $r$, identify every class $\mathcal{E}_{\mathbf{\bar w_k}, m}^r$ by a minimal representation.}\\
\hline
\textbf{Step 1} & Generate all $J:=(m!)^n$ profiles $\mathbf{P}^{1},\dots,\mathbf{P}^{J}\in\mathcal{P}(A)^n$ for $A:=\{a_1,\ldots, a_m\}$. Set $\mathcal{F}:=\varnothing$. \\
\textbf{Step 2} & For every $\mathbf{P}^{j} \in \mathcal{P}(A)^n$ and every $a_i \in A$, check if there is any weight vector $\mathbf{w}%=(w_1,...,w_n) 
\in \mathbb{N}_0^n$\ \,s.t.\ $r|\mathbf{w}(\mathbf{P}^{j})=a_i$ by testing feasibility of the implied constraints (cf.\ Section~\ref{subsec:illustration}). 
If yes, then append $(i,j)$ to $\mathcal{F}$.\\
\textbf{Step 3}  & Loop over $j$ from $1$ to $J$.\\                       
\textbf{Step 3a} & If $j=1$, then set $\mathcal{C}_1:=\left\{1\le i\le m \ |\  (i,j)\in\mathcal{F}\right\}$. \\                       
\textbf{Step 3b} & If $j\ge 2$, then set $\mathcal{C}_j:=\varnothing$ and loop over all $(p_1,\dots,p_{j-1})\in\mathcal{C}_{j-1}$ and all
$p_j\in\{1,\dots,m\}$ with $(p_j,j)\in\mathcal{F}$. If (\ref{eq:ilp_min_sum}) has a solution for the restriction
to the profiles $\mathbf{P^{1}},\dots,\mathbf{P^{j}}$ 
with prescribed winners $\rho(\mathbf{P^{i}})=a_{p_i}$ for $1\le i\le j$,
then append $(p_1,\dots,p_j)$ to $\mathcal{C}_p$. \\
\textbf{Step 4}  & Loop over the elements $(p_1, \ldots, p_j,\ldots, p_J)\in \mathcal{C}_J$ and output minimal weights $\mathbf{\bar w}$ such that $r|\mathbf{\bar w}\equiv\rho$ 
with $\rho(\mathbf{P^j})=p_j$ by solving (\ref{eq:ilp_min_sum}). \\
\hline
\end{tabular}\caption{Determining the classes of $r$-weighted committees for given $n$ and $m$ \label{tab:BC_algorithm}
} 
\end{table}

The algorithm can still require impractical memory size and running time.  
The main alternative then is to heuristically loop over different weight distributions and check if they are structurally distinct from those already known. Namely,  
start with $w_\Sigma:=0$ and an empty list $\hat{\mathcal{W}}$ of weight vectors; increase the sum of weights $w_\Sigma$ in steps of $1$; generate the set 
%\begin{equation}\label{eq:heuristic}
$\mathcal{W}_{w_\Sigma}:=\big\{\mathbf{w} \in \mathbb{N}_0^n \ \big|\  w_1 \geq \dots \geq w_n \text{ and } 
w_1+\dots+w_n=w_\Sigma\big\}
$
and loop over all $\mathbf{w}\in \mathcal{W}_{w_\Sigma}$.
The respective weight vector $\mathbf{w}$ is appended to $\hat{\mathcal{W}}$ if for every $\mathbf{w'}\in \hat{\mathcal{W}}$ we have $r|\mathbf{w}(\mathbf{P})\neq r|\mathbf{w'}(\mathbf{P})$ for at least one 
$\mathbf{P}\in\mathcal{P}(A)^n$. 
The set $\hat{\mathcal{W}}$ then contains a growing list of minimal weight vectors that induce different mappings from preference profiles to winners and hence correspond to structurally distinct committee games $[N,A,r,\mathbf{w}]$. 
This method has the advantage of not requiring a weightedness test, such as (ILP). However,
search needs to be stopped manually and just produces
a lower bound on the actual number of classes.\footnote{Upper bounds exist for weight sums that guarantee coverage of \emph{all} equivalence classes (cf.\ \citeNP[Thm.\ 9.3.2.1]{Muroga:1971}). 
The bounds are too large to be practical, however.} 

\section{Number and geometry of weighted committee games}\label{sec:results}

\subsection{Number of antiplurality, Borda, Copeland, and plurality games}

\renewcommand{\arraystretch}{1.2}

\begin{table}%[htbp] 
	\begin{center}
		\begin{tabular}{||c|c|c|c|c||}
			\hline \hline
			\diagbox[innerleftsep=0cm,innerrightsep=0cm]{\ $n,m$}{\\ $r$\ } & Antiplurality &        Borda        &       Copeland        &      Plurality       \\
			\hline\hline
			{3,2} &    \multicolumn{1}{l}{}&      \multicolumn{2}{c}{4}        & \multicolumn{1}{c||}{} \\
			{4,2} &    \multicolumn{1}{l}{}&      \multicolumn{2}{c}{9}        & \multicolumn{1}{c||}{} \\
			{5,2} &    \multicolumn{1}{l}{}&     \multicolumn{2}{c}{27}        & \multicolumn{1}{c||}{} \\
			{6,2} &    \multicolumn{1}{l}{}&     \multicolumn{2}{c}{138}       & \multicolumn{1}{c||}{} \\
			{7,2} &    \multicolumn{1}{l}{}&   \multicolumn{2}{c}{1\,663}      & \multicolumn{1}{c||}{} \\
			{8,2} &    \multicolumn{1}{l}{}&   \multicolumn{2}{c}{63\,764}     & \multicolumn{1}{c||}{} \\
			{9,2} &    \multicolumn{1}{l}{}& \multicolumn{2}{c}{9\,425\,479}   & \multicolumn{1}{c||}{} \\
			\hline                                                                                 
			$3,3$ &       5       &         51          &           4           &          6           \\
			$3,4$ &       3       &         505         &           4           &          6           \\
			$3,5$ &       3       &   {$\geq 2\,251$}   &           4           &          6           \\
			\hline
			$4,3$ &      19       &       5\,255        &           9           &          34          \\
			$4,4$ &       7       &  {$\gg 635\,622$}   &           9           &          36          \\
			$4,5$ &       4       &  {$\gg 635\,622$}   &           9           &          36          \\
			\hline
			$5,3$ &      263      & {$\gg 1\,153\,448$} &          27           &         852          \\
			\hline
			$6,3$ &      $\geq 33\,583$     & {$\gg 1\,153\,448$} &          138           &         {$\gg 147\,984$}          \\
			\hline \hline
		\end{tabular}
	\end{center}
	\caption{\small Numbers of distinct weighted committee games 
	}
	\label{table:nr_classes}
\end{table}

A combination of our analytical observations and computational means permits identification of all structurally distinct $r$-weighted committee games with $r\in \{r^A, r^B, r^C, r^P\}$ for small $n$ and $m$. 
% Klassen für m=2:
% Plurality, Copeland, Antiplurality: Für n=2 und m>=2 immer 2 Klassen
% Borda: m=3: 4, m=4: 8, m=5: 12, m=6: 20, m=7: 24 
Table~\ref{table:nr_classes} summarizes our findings; 
figures do not include $\mathcal{E}_{\mathbf{0}, m}$.\footnote{$\mathbf{w}^0=\mathbf{0}$ always forms its own equivalence class:
	consider the unanimous profile $\mathbf{P}=(P,\ldots, P)\in\mathcal{P}(A)^n$ with  $a_2 P a_3 P \dots P a_m P a_1$. Then $r|\mathbf{0}(\mathbf{P})=a_1$ but $r|\mathbf{w}(\mathbf{P})=a_2$ for any $\mathbf{w}\neq \mathbf{0}$.} 
The branch-and-cut approach required excessive memory for Borda committees when $m>4$ or $n=m\ge 4$.\footnote{We used 128 GB RAM and eight 3.0 GHz cores. Several instances ran for %several 
	more than six months.}
We write ``$\ge \ldots$'' if the heuristic appended no new games to set $\hat{\mathcal{W}}$ for long enough to support the conjecture that the reported bound equals the exact number of games; we write ``$\gg \ldots$'' otherwise. 

Figures for $m=2$ and $n\le 6$ have been obtained in the literature before; the others are, to our knowledge, new.
When less than 150 equivalence classes of games exist, we report minimal representations in Appendix~B.
Our list for $m=2$ nests the weighted voting games with 50\%-majority threshold reported by \citeN{Krohn/Sudholter:1995} and \citeN{Brams/Fishburn:1996}; plurality committees with $m=3$ nest the subset of tie-free games identified by \shortciteN{Chua/Ueng/Huang:2002} for $n=3$, 4.

%%% HINWEIS: Die Zahlen bei Brams/Fishburn (1996) unterscheiden sich von unseren. Dies liegt an: 1): BF schließen das Diktaturspiel aus. 
%%% 2) BF erlauben kein Nullgewicht. Das führt dazu, dass B/F Spiele "verlieren". Bei 3 Spielern fehlt ihnen z.B. das Spiel 
%%% [1 1 0], da es für n=3 und q>0.5 unmöglich ist, einen Dummy Spieler zu haben, wenn gleichzeitig jeder Spieler mindestens eine 
%%% Stimme haben muss. Bei n=4 gilt das gleiche für den Fall von 2 Dummy Spielern. Daher fehlt bei BF das Spiel [1 1 0 0]. Der Fall mit 
%%% nur einem Dummy Spielern bei n=4 wird bei BF durch das Spiel [2 2 2 1] abgedeckt. Wird aber ein Nullgewicht zugelassen, gibt 
%%% es zwei Möglichkeiten einen Dummy Spieler zu haben: [2 1 1 0] und [1 1 1 0]. Nehmen wir aus unserer Liste jedoch zunächst das 
%%% Diktaturspiel und alle Spiele, in denen ein Spieler ein Nullgewicht hat raus und fügen anschließend die Spiele, die bei BF einen 
%%% Dummy Spieler induzieren wieder hinzu, kommen wir - bis auf n=6 - auf die gleichen Zahlen wie BF. Bei n=6 erhalten wir nach dieser 
%%% "Methode" 117 Spiele, BF jedoch nur 116. Ein Vergleich der zugehörigen SSIs liefert leider auch keine Klarheit: Unsere Spiele 
%%% [9 7 5 4 3 2], [5 4 4 3 3 1] und [7 6 5 4 4 2] (bzw. deren SSI) kommen bei BF nicht vor. Dagegen kommen die Spiele [4 3 2 2 1 1] 
%%% und [4 3 3 2 2 1] (bzw. deren SSI) bei BF doppelt vor. Ersteres taucht in Form von [29 22 20 19 5 5] (Nr. 41) und [30 21 19 19 9 2] 
%%% (Nr. 58) auf. Letzteres in Form von [27 22 20 15 9 7] (Nr. 90) und [27 22 20 15 12 4] (Nr. 92).   

\subsection{Geometry of committee games with \emph{n}=3}

\begin{figure}
	{\small (a)  \hspace{0.495\textwidth}(b)}  \\ 
	\includegraphics[width=0.49\textwidth]{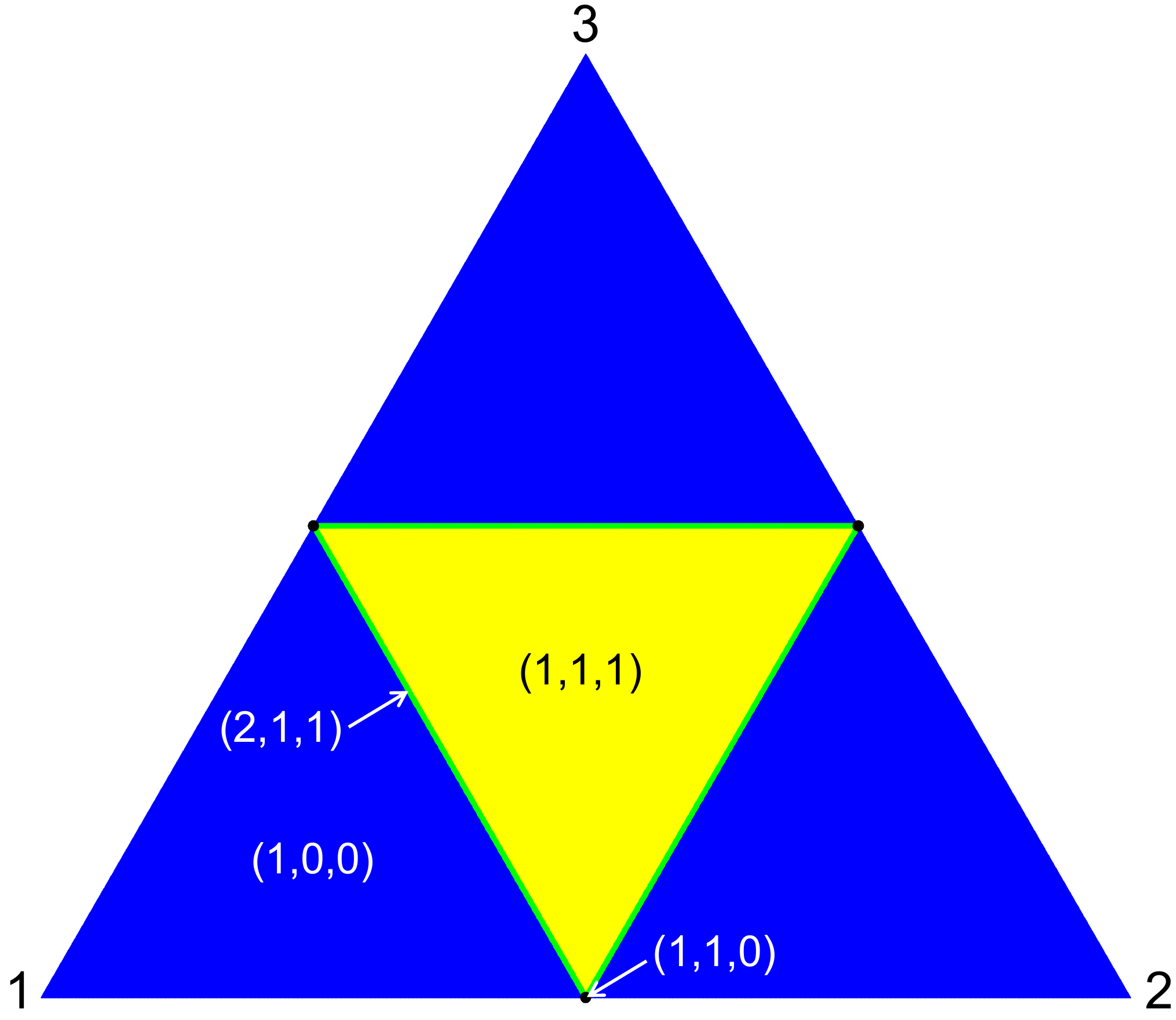}
	\includegraphics[width=0.49\textwidth]{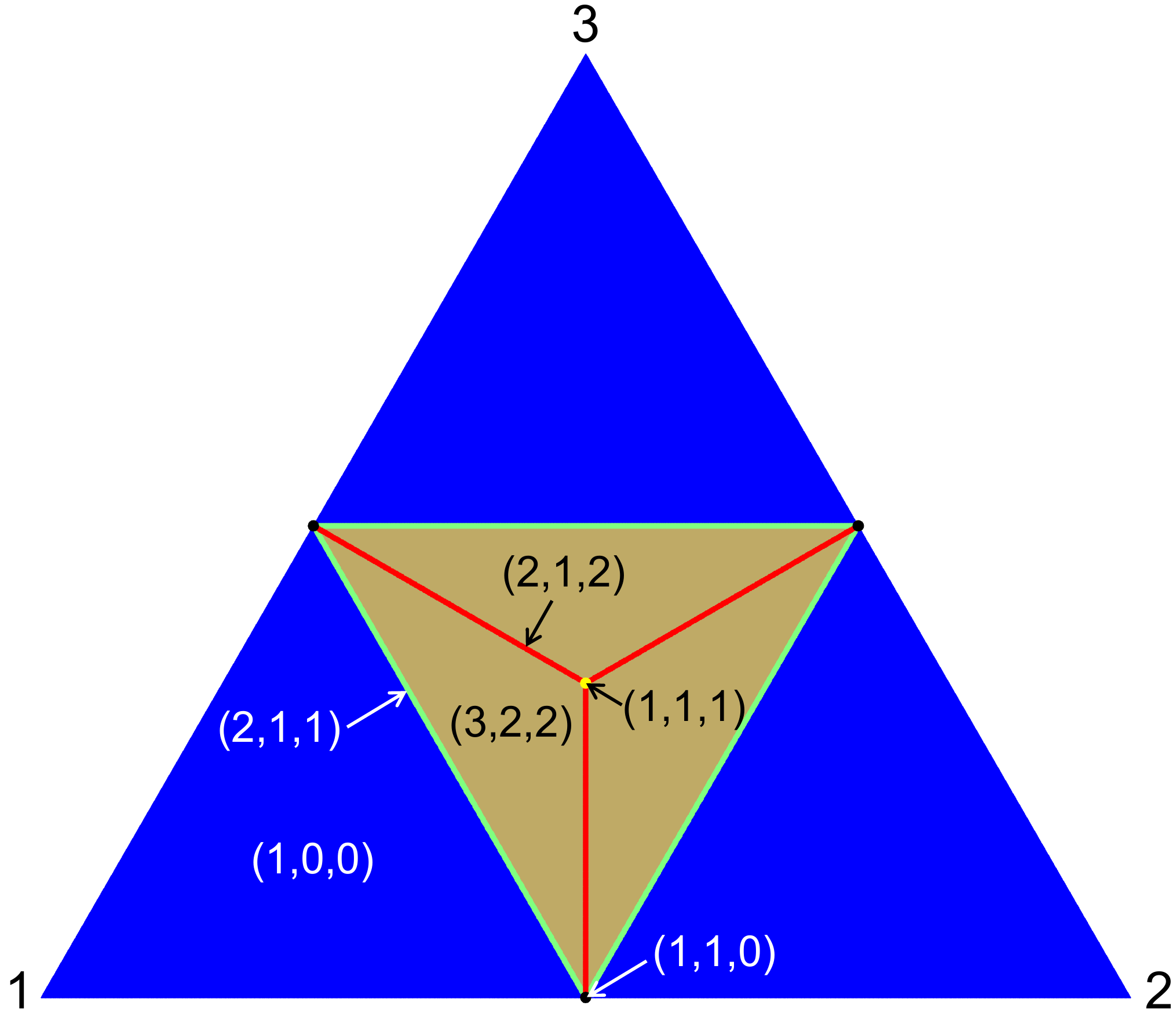}
	\caption{\small The four Copeland and six plurality equivalence classes}
	\label{fig:EQ_Copeland_Plurality}
\end{figure}

\begin{figure}
	{\small (a) $m=3$ \hspace{0.495\textwidth}(b) $m\ge 4$}  \\ 
  \includegraphics[width=0.49\textwidth]{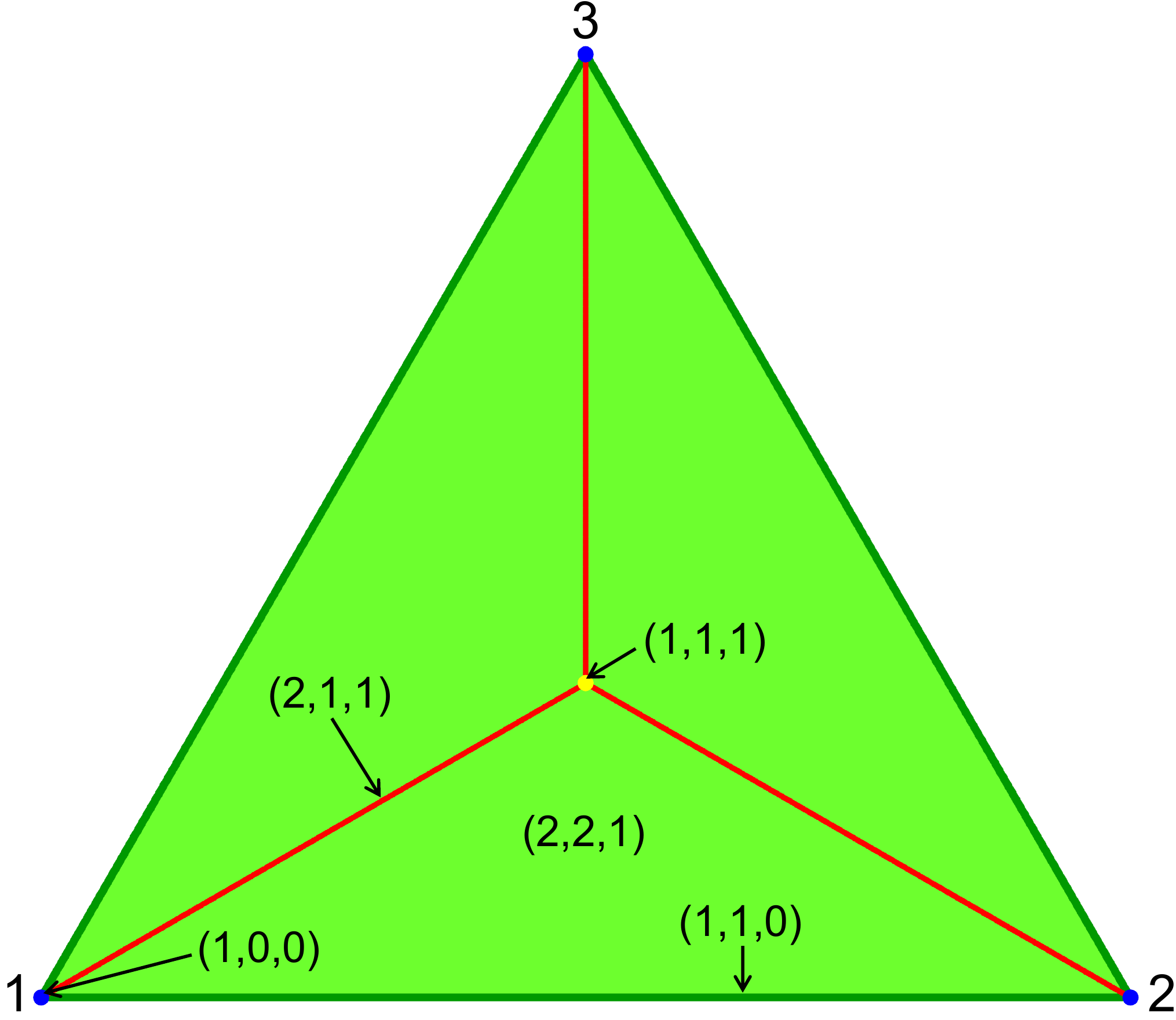}
  \includegraphics[width=0.49\textwidth]{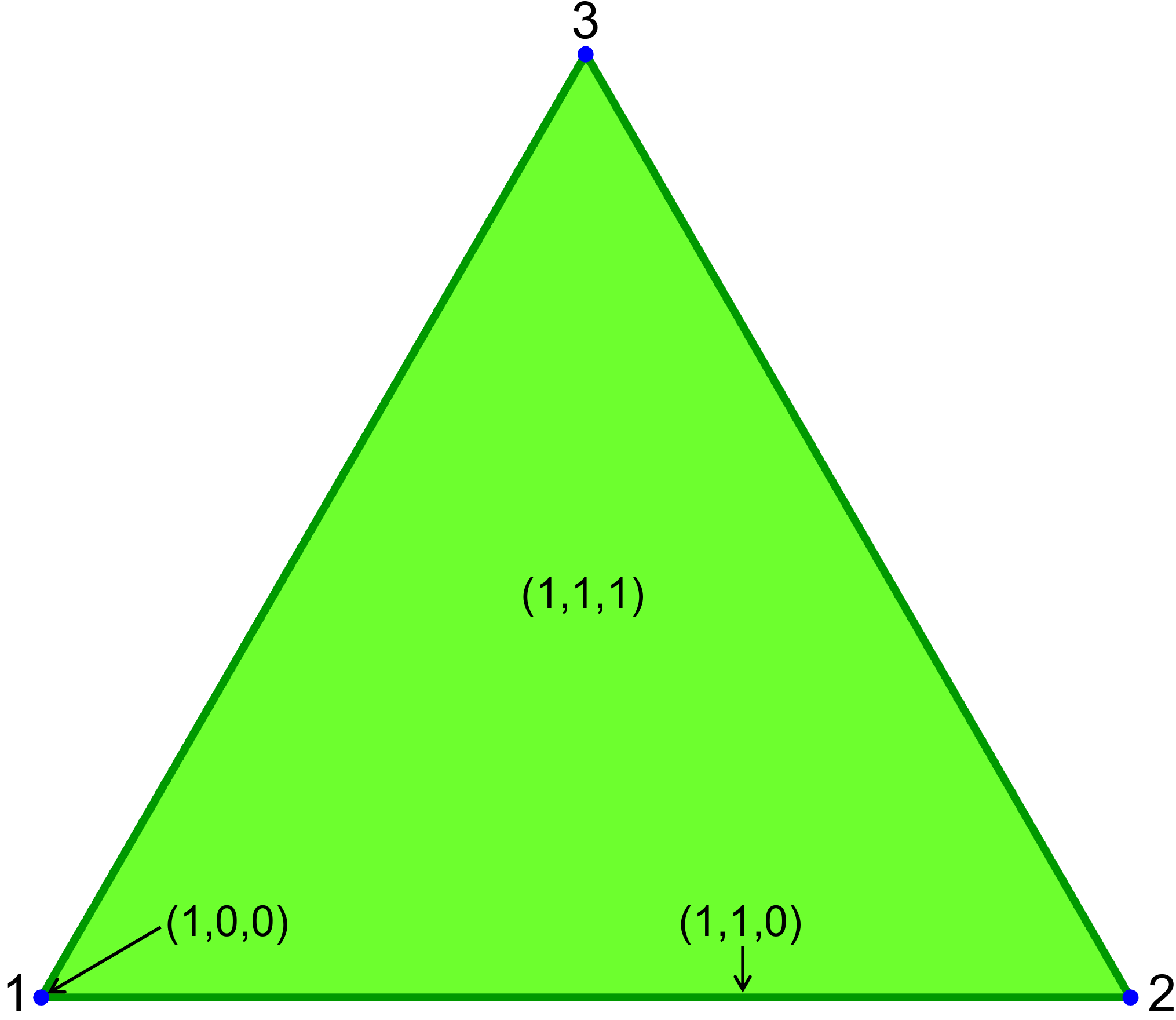}
      \caption{\small The five or three antiplurality equivalence classes}
\label{fig:EQ_Antiplurality}
\end{figure}

The case of three players allows to convey a geometric idea of how robust a given weighted committee is to small changes in weights. 
Our illustrations echo those by Saari (\citeyearNP{Saari:1995}, \citeyearNP{Saari:2001:chaotic}): his eponymous triangles concern $m=3$ alternatives and arbitrary numbers $n$ of individual voters. 
They illuminate how collective rankings vary with the applicable voting procedure for fixed preferences $\mathbf{P}$. 
We, by contrast, assume $n=3$ voter blocs, evaluate all preferences, and the number~$m$ of alternatives may vary. 
We use the standard projection of the 3-dimensional unit simplex of relative weights to the plane. 
Points of identical color correspond to structurally equivalent weight distributions, i.e., they induce isomorphic committee games for the investigated
voting rule $r$. 
When equivalence classes are line segments or single points, we have manually enlarged them in Figures~\ref{fig:EQ_Copeland_Plurality} and \ref{fig:EQ_Antiplurality} to improve visibility.

Figure~\ref{fig:EQ_Copeland_Plurality}(a) shows all Copeland committees with three players and their minimal representations. 
Very dissimilar weight distributions like $(33, 33, 33)$ and $(49, 49, 1)$ induce the same Copeland winners. 
Figure~\ref{fig:EQ_Copeland_Plurality}(b) illustrates that plurality rule is more sensitive to weight perturbations than Copeland rule, at least for non-dictatorial configurations. This is more pronounced the more voter blocs are involved: 
there are about four and 32 times more distinct committees with plurality than Copeland rule for $n=4$ and 5 (Table~\ref{table:nr_classes}); we conjecture the factor exceeds 1\,000 for $n=6$.

\begin{figure}
	{\small \mbox{ }  \hspace{0.325\textwidth} (a) $m=3$} \\
	\mbox{ }  \hspace{0.23\textwidth}	\includegraphics[width=0.49\textwidth]{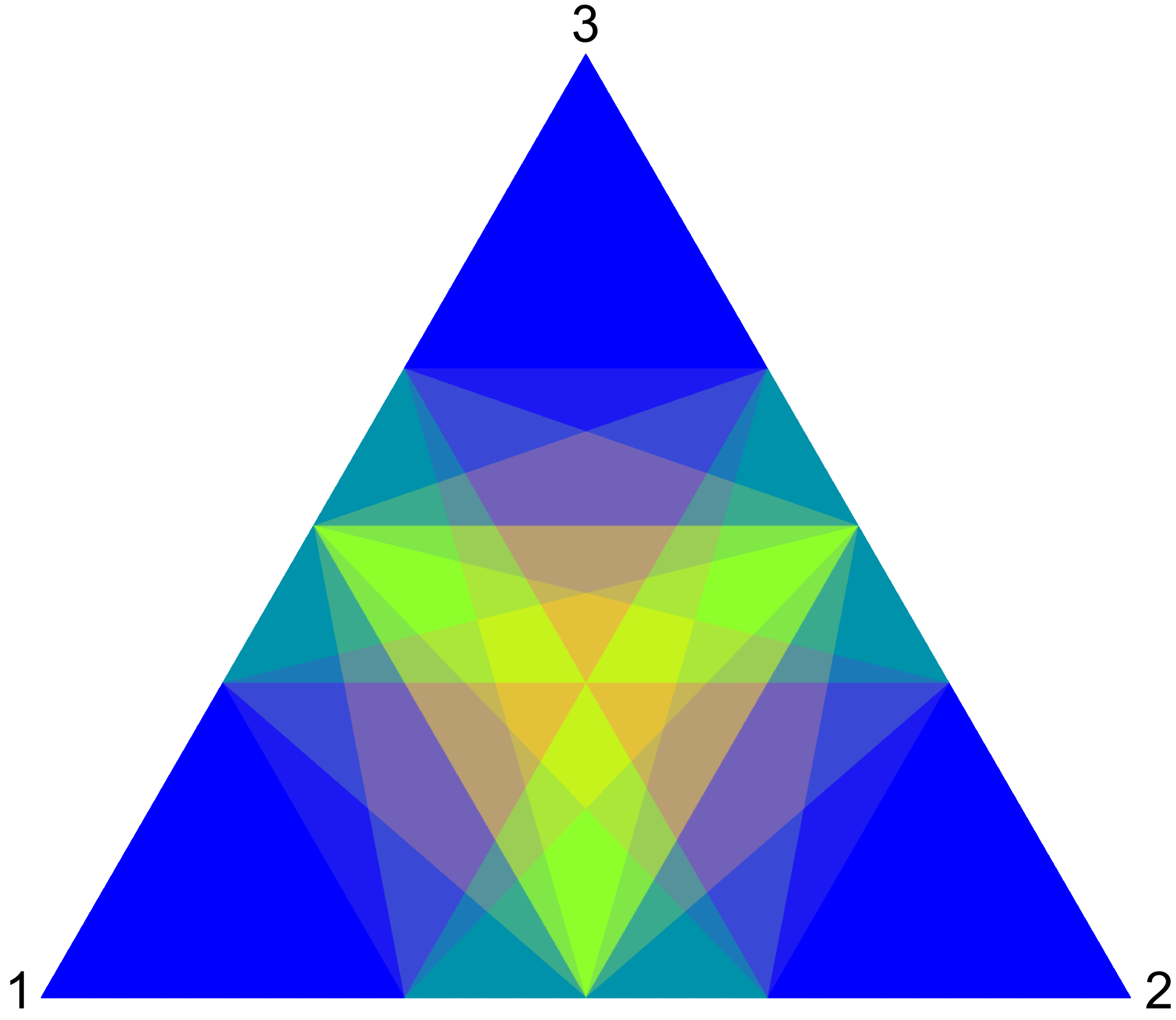} \\
	%	\mbox{ } \\
	{\small \mbox{ }  \hspace{0.08\textwidth} (b) $m=4$ \hspace{0.39\textwidth}(c) $m= 5$}  \\ 
	\includegraphics[width=0.49\textwidth]{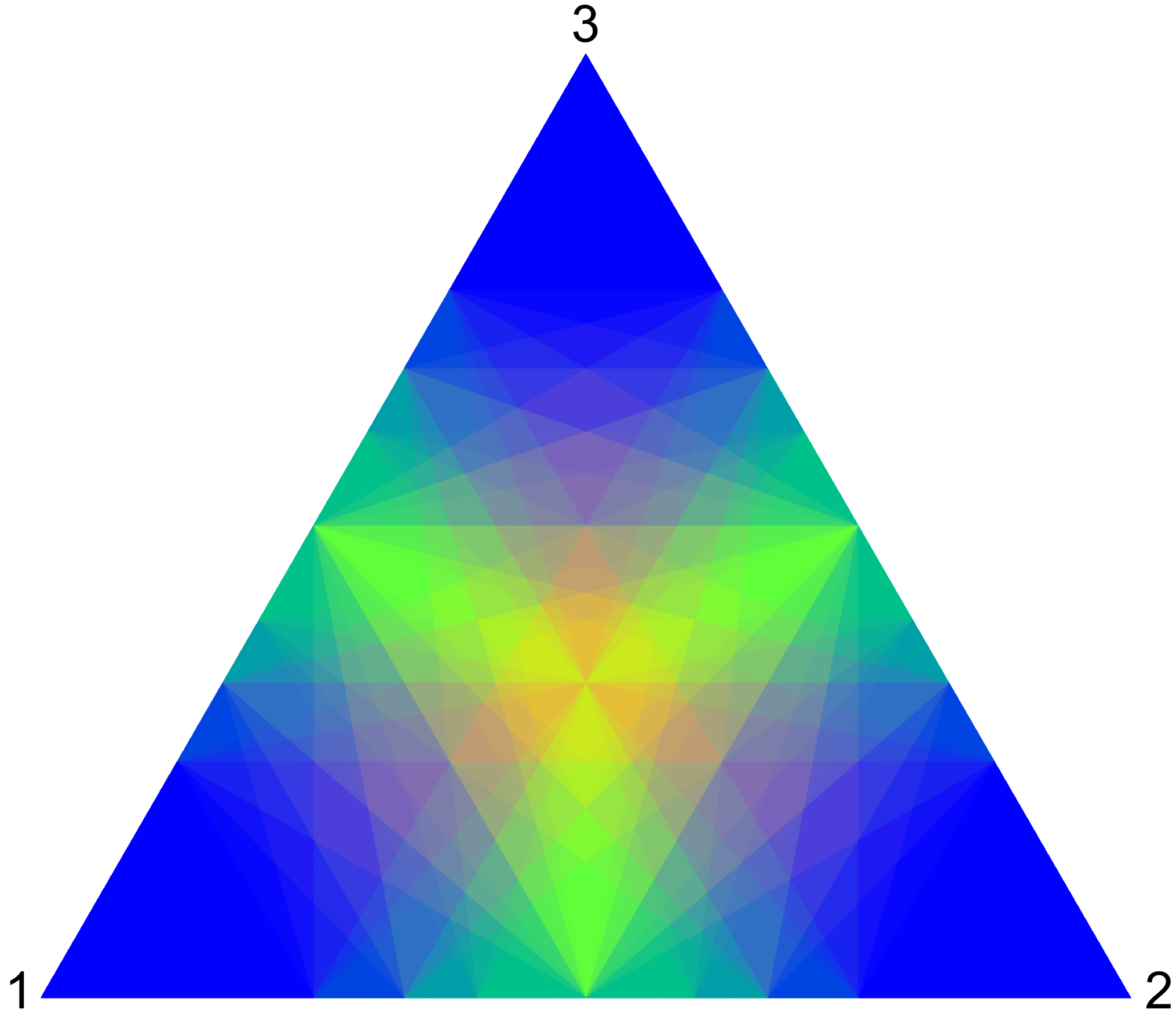}
	\includegraphics[width=0.49\textwidth]{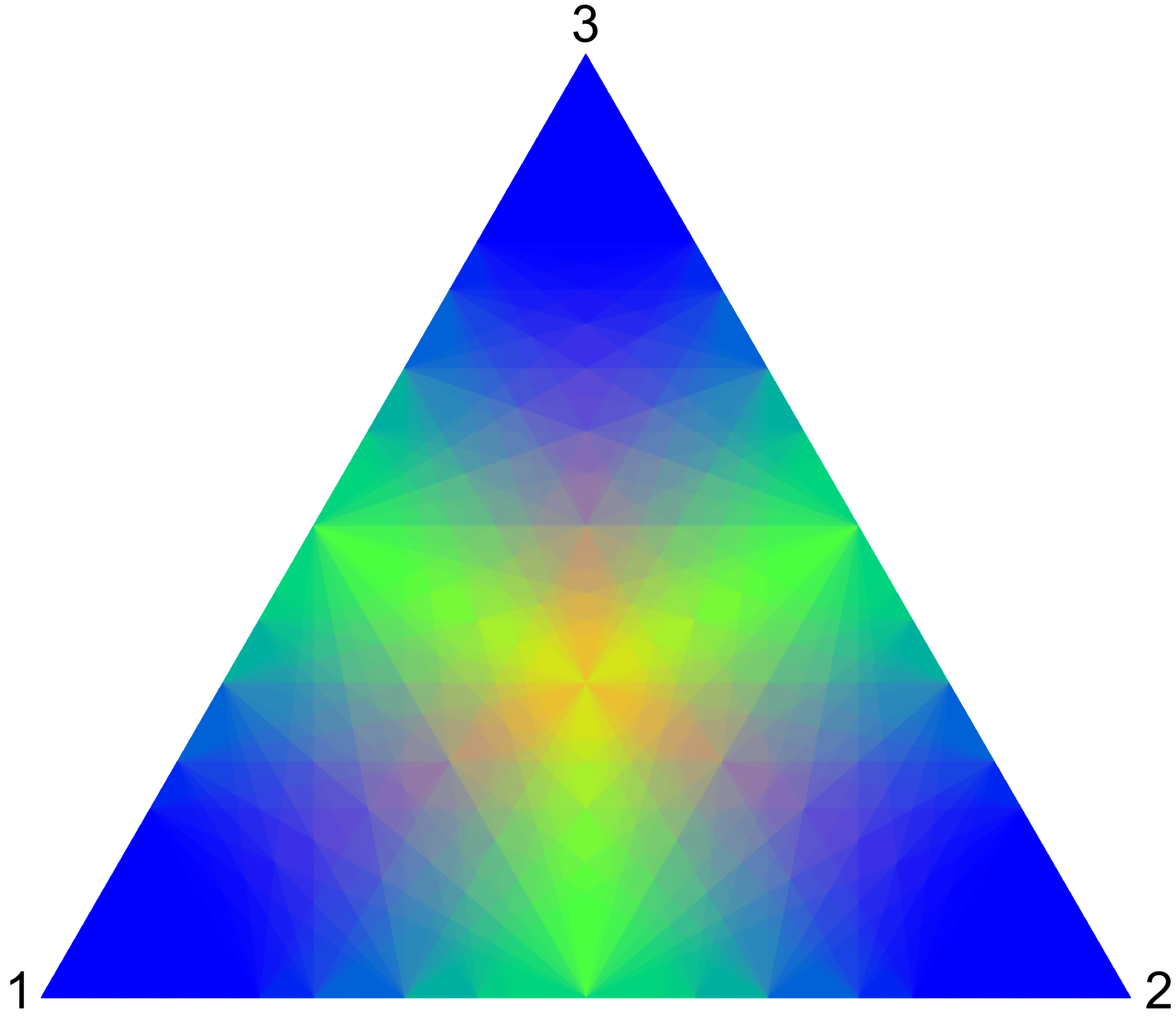} \\
	
	\caption{\small The 51, 505, and $\ge\! 2251$ Borda equivalence classes}
	\label{fig:EQ_Borda}
\end{figure}

Antiplurality rule (Figure~\ref{fig:EQ_Antiplurality}), characterized recently with new axioms by \citeN{Kurihara:2018}, holds an intermediate ground in terms of the sensitivity to $\mathbf{w}$. 
The most  scope for changes in the distribution of voting rights to induce different decisions
comes with Borda rule, as illustrated by Figure~\ref{fig:EQ_Borda}. 
(Note that Figure~\ref{fig:EQ_Copeland_Plurality}(a) captures the case of $m=2$ for $r^A$, $r^B$, and $r^P$, too.) 
Whether such sensitivity is (un)desir\-able from an institutional perspective will depend on context and objectives. 
Higher sensitivity can give bigger incentives for political parties to campaign or private investment in voting stock. 
However, this needs to be weighed against other properties of a voting method such as consistency with respect to subgroup decisions \cite{Young:1975}, informational requirements, or complexity of strategic manipulation.

\section{Concluding remarks}\label{sec:conclusion}
Equivalence of seemingly different distributions of seats, quotas, voting stock, etc.\ 
depends highly on whether decisions involve two, three, or more alternatives.
Weight distributions such as $(6,5,2)$, $(1,1,1)$, or $(48\%, 24\%, 28\%)$ 
induce the same binary majority choices and Copeland decisions but lead to non-equivalent mappings from preferences to selected alternatives in other cases. 
The respective potential for weight differences to matter has been formalized and compared here. 
The investigated equivalences and their voting power implications (see \shortciteNP{Kurz/Mayer/Napel:2019}) could be of interest not only for multicandidate voting in corporations, councils, or parliaments; e.g., acceptable polling error should be larger in the middle of an equivalence class of an election rule than on its boundary.

For the IMF's Executive Board, we have checked that the 2016 reform of drawing rights can have consequences for who becomes the next IMF~Managing Director. Election winners differ between reformed and unreformed weights for about 5\% of all $6^{24} $ %\approx 4.7\cdot 10^{18}$ 
conceivable strict preference configurations over three shortlisted candidates.
It matters for almost 15\% of profiles whether Copeland or plurality rule is applied.

There is ample choice for extending above analysis since the list of single-winner voting procedures used by committees is long (see, e.g., \citeNP{Laslier:2012}). 
We have tentatively considered scoring rules with arbitrary 
$\mathbf{s}=(1, s_2, 0)\in \mathbb{Q}^3$ for $n=m=3$, too.  
The numbers of distinct committees are M-shaped: they increase from 6 for plurality to more than 160 for $s_2=0.25$,  %163
fall to 51 Borda committees for $s_2=0.5$, increase again to at least 229 for $s_2=0.9$ and then drop to 5 antiplurality committees. % for $s_2=1$. %\footnote{Illustrations of their geometry are available upon request. Some reminded us of %are almost reminiscences of paintings, e.g., by Bauhaus artists Paul Klee and Johannes Itten.
%} 
%

Implications of different weights are, obviously, just one aspect of preference aggregation by voting among others (cf.\ \citeNP{Nurmi:1987}). 
Group decisions may often be purely binary affairs and the fast-growing dimensionality of collective choice makes it challenging to evaluate 
weighted voting for $m\ge 3$.  
Table~\ref{table:nr_classes} and Figures~\ref{fig:EQ_Copeland_Plurality}--\ref{fig:EQ_Borda} document how non-trivial the links between 
weight and choice differences %distributions and different choices 
can be. They seem relevant enough, however, to be studied %not just for the % 
beyond the 
binary case.

%%%%%%%%%%%%%%%%%%%%%
%%%%%%%%%%%%%%%%%%%%%
%%%%%%%%%%%%%%%%%%%%% REFERENCES
%%%%%%%%%%%%%%%%%%%%%
%%%%%%%%%%%%%%%%%%%%%
%\newpage 
\small

\phantomsection \label{References}
\addcontentsline{toc}{section}{References}

\setstretch{1.18} 

\setlength{\labelsep}{-0.2cm}

%\bibliographystyle{chicago}
%\bibliography{KMN-Literatur}

\begin{thebibliography}{}
	
	\bibitem[\protect\citeauthoryear{Aleskerov and Kurbanov}{Aleskerov and
		Kurbanov}{1999}]{Aleskerov/Kurbanov:1999}
	Aleskerov, F. and E.~Kurbanov (1999).
	\newblock Degree of manipulability of social choice procedures.
	\newblock In A.~Alkan, C.~D. Aliprantis, and N.~C. Yannelis (Eds.), {\em
		Current Trends in Economics}, pp.\  13--27. Berlin: Springer.
	
	\bibitem[\protect\citeauthoryear{Amer, Carreras, and Mag{\~{a}}na}{Amer
		et~al.}{1998}]{Amer/Carreras/Magana:1998b}
	Amer, R., F.~Carreras, and A.~Mag{\~{a}}na (1998).
	\newblock Extension of values to games with multiple alternatives.
	\newblock {\em Annals of Operations Research\/}~{\em 84\/}(0), 63--78.
	
	\bibitem[\protect\citeauthoryear{Bolger}{Bolger}{1983}]{Bolger:1983}
	Bolger, E.~M. (1983).
	\newblock The {Banzhaf} index for multicandidate presidential elections.
	\newblock {\em SIAM Journal on Algebraic and Discrete Methods\/}~{\em 4\/}(4),
	442--458.
	
	\bibitem[\protect\citeauthoryear{Bouton}{Bouton}{2013}]{Bouton:2013}
	Bouton, L. (2013).
	\newblock A theory of strategic voting in runoff elections.
	\newblock {\em American Economic Review\/}~{\em 103\/}(4), 1248--1288.
	
	\bibitem[\protect\citeauthoryear{Brams and Fishburn}{Brams and
		Fishburn}{1996}]{Brams/Fishburn:1996}
	Brams, S.~J. and P.~C. Fishburn (1996).
	\newblock Minimal winning coalitions in weighted-majority voting games.
	\newblock {\em Social Choice and Welfare\/}~{\em 13\/}(4), 397--417.
	
	\bibitem[\protect\citeauthoryear{Cheung and Ng}{Cheung and
		Ng}{2014}]{Cheung/Ng:2014}
	Cheung, W.-S. and T.-W. Ng (2014).
	\newblock A three-dimensional voting system in Hong Kong.
	\newblock {\em European Journal of Operational Research\/}~{\em 236\/}(1),
	292--297.
	
	\bibitem[\protect\citeauthoryear{Chua, Ueng, and Huang}{Chua
		et~al.}{2002}]{Chua/Ueng/Huang:2002}
	Chua, V. C.~H., C.~H. Ueng, and H.~C. Huang (2002).
	\newblock A method for evaluating the behavior of power indices in weighted
	plurality games.
	\newblock {\em Social Choice and Welfare\/}~{\em 19\/}(3), 665--680.
	
	\bibitem[\protect\citeauthoryear{Felsenthal and Machover}{Felsenthal and
		Machover}{1997}]{Felsenthal/Machover:1997}
	Felsenthal, D.~S. and M.~Machover (1997).
	\newblock Ternary voting games.
	\newblock {\em International Journal of Game Theory\/}~{\em 26\/}(3), 335--351.
	
	\bibitem[\protect\citeauthoryear{Felsenthal and Nurmi}{Felsenthal and
		Nurmi}{2017}]{Felsenthal/Nurmi:2017}
	Felsenthal, D.~S. and H.~Nurmi (2017).
	\newblock {\em Monotonicity Failures Afflicting Procedures for Electing a
		Single Candidate}.
	\newblock Cham: Springer.
	
	\bibitem[\protect\citeauthoryear{Freixas, Freixas, and Kurz}{Freixas
		et~al.}{2017}]{Freixas/Freixas/Kurz:2017}
	Freixas, J., M.~Freixas, and S.~Kurz (2017).
	\newblock On the characterization of weighted simple games.
	\newblock {\em Theory and Decision\/}~{\em 83\/}(4), 469--498.

	\bibitem[\protect\citeauthoryear{Freixas and Kaniovski}{Freixas and
		Kaniovski}{2014}]{Freixas/Kaniovski:2014}
	Freixas, J. and S.~Kaniovski (2014).
	\newblock The minimum sum representation as an index of voting power.
	\newblock {\em European Journal of Operational Research\/}~{\em 233\/}(3),
	739--748.
	
	\bibitem[\protect\citeauthoryear{Freixas and Zwicker}{Freixas and
		Zwicker}{2003}]{Freixas/Zwicker:2003}
	Freixas, J. and W.~S. Zwicker (2003).
	\newblock Weighted voting, abstention, and multiple levels of approval.
	\newblock {\em Social Choice and Welfare\/}~{\em 21\/}(3), 399--431.
	
	\bibitem[\protect\citeauthoryear{Henriet}{Henriet}{1985}]{Henriet:1985}
	Henriet, D. (1985).
	\newblock The {Copeland} choice function: an axiomatic characterization.
	\newblock {\em Social Choice and Welfare\/}~{\em 2\/}(1), 49--63.
	
	\bibitem[\protect\citeauthoryear{Houy and Zwicker}{Houy and
		Zwicker}{2014}]{Houy/Zwicker:2014}
	Houy, N. and W.~S. Zwicker (2014).
	\newblock The geometry of voting power: weighted voting and hyper-ellipsoids.
	\newblock {\em Games and Economic Behavior\/}~{\em 84}, 7--16.
	
	\bibitem[\protect\citeauthoryear{Hsiao and Raghavan}{Hsiao and
		Raghavan}{1993}]{Hsiao/Raghavan:1993}
	Hsiao, C.-R. and T.~E.~S. Raghavan (1993).
	\newblock {Shapley} value for multichoice cooperative games, {I}.
	\newblock {\em Games and Economic Behavior\/}~{\em 5\/}(2), 240--256.
	
	\bibitem[\protect\citeauthoryear{Krohn and Sudh\"{o}lter}{Krohn and
		Sudh\"{o}lter}{1995}]{Krohn/Sudholter:1995}
	Krohn, I. and P.~Sudh\"{o}lter (1995).
	\newblock Directed and weighted majority games.
	\newblock {\em Mathematical Methods of Operations Research\/}~{\em 42\/}(2),
	189--216.
	
	\bibitem[\protect\citeauthoryear{Kurihara}{Kurihara}{2018}]{Kurihara:2018}
	Kurihara, T. (2018).
	\newblock A simple characterization of the anti-plurality rule.
	\newblock {\em Economics Letters\/}~{\em 168}, 110--111.
	
	\bibitem[\protect\citeauthoryear{Kurz}{Kurz}{2012a}]{Kurz:2012:representation}
	Kurz, S. (2012a).
	\newblock On minimum sum representations for weighted voting games.
	\newblock {\em Annals of Operations Research\/}~{\em 196\/}(1), 361--369.
	
	\bibitem[\protect\citeauthoryear{Kurz}{Kurz}{2012b}]{Kurz:2012}
	Kurz, S. (2012b).
	\newblock On the inverse power index problem.
	\newblock {\em Optimization\/}~{\em 61\/}(8), 989--1011.
	
	\bibitem[\protect\citeauthoryear{Kurz, Mayer, and Napel}{Kurz
		et~al.}{2019}]{Kurz/Mayer/Napel:2019}
	Kurz, S., A.~Mayer, and S.~Napel (2019).
	\newblock Influence in weighted committees.
	\newblock Mimeo, University of Bayreuth.	
	
	\bibitem[\protect\citeauthoryear{Kurz and Napel}{Kurz and
		Napel}{2016}]{Kurz/Napel:2016}
	Kurz, S. and S.~Napel (2016).
	\newblock Dimension of the {Lisbon} voting rules in the {EU Council}: A
	challenge and new world record.
	\newblock {\em Optimization Letters\/}~{\em 10\/}(6), 1245--1256.
	
	\bibitem[\protect\citeauthoryear{Laruelle and Valenciano}{Laruelle and
		Valenciano}{2012}]{Laruelle/Valenciano:2012}
	Laruelle, A. and F.~Valenciano (2012).
	\newblock Quaternary dichotomous voting rules.
	\newblock {\em Social Choice and Welfare\/}~{\em 38\/}(3), 431--454.
	
	\bibitem[\protect\citeauthoryear{Laslier}{Laslier}{2012}]{Laslier:2012}
	Laslier, J.-F. (2012).
	\newblock And the loser is \ldots plurality voting.
	\newblock In D.~S. Felsenthal and M.~Machover (Eds.), {\em Electoral Systems:
		Paradoxes, Assumptions, and Procedures}, pp.\  327--351. Berlin: Springer.
	
	\bibitem[\protect\citeauthoryear{Monroy and Fern\'{a}ndez}{Monroy and
		Fern\'{a}ndez}{2009}]{Monroy/Fernandez:2009}
	Monroy, L. and F.~R. Fern\'{a}ndez (2009).
	\newblock A general model for voting systems with multiple alternatives.
	\newblock {\em Applied Mathematics and Computation\/}~{\em 215\/}(4),
	1537--1547.
	
	\bibitem[\protect\citeauthoryear{Monroy and Fern\'{a}ndez}{Monroy and
		Fern\'{a}ndez}{2011}]{Monroy/Fernandez:2011}
	Monroy, L. and F.~R. Fern\'{a}ndez (2011).
	\newblock The {Shapley-Shubik} index for multi-criteria simple games.
	\newblock {\em European Journal of Operational Research\/}~{\em 209\/}(2),
	122--128.
	
	\bibitem[\protect\citeauthoryear{Moulin}{Moulin}{1988}]{Moulin:1988}
	Moulin, H. (1988).
	\newblock Condorcet's principle implies the no show paradox.
	\newblock {\em Journal of Economic Theory\/}~{\em 45\/}(1), 53--64.
	
	\bibitem[\protect\citeauthoryear{Muroga}{Muroga}{1971}]{Muroga:1971}
	Muroga, S. (1971).
	\newblock {\em Threshold Logic and its Applications}.
	\newblock New York, NY: Wiley.
	
	\bibitem[\protect\citeauthoryear{Myerson and Weber}{Myerson and
		Weber}{1993}]{Myerson/Weber:1993}
	Myerson, R.~B. and R.~J. Weber (1993).
	\newblock A theory of voting equilibria.
	\newblock {\em American Political Science Review\/}~{\em 87\/}(1), 102--114.
	
	\bibitem[\protect\citeauthoryear{Nurmi}{Nurmi}{1987}]{Nurmi:1987}
	Nurmi, H. (1987).
	\newblock {\em Comparing Voting Systems}.
	\newblock Dordrecht: Kluwer.
	
	\bibitem[\protect\citeauthoryear{Saari}{Saari}{1995}]{Saari:1995}
	Saari, D.~G. (1995).
	\newblock {\em Basic Geometry of Voting}.
	\newblock Berlin: Springer.
	
	\bibitem[\protect\citeauthoryear{Saari}{Saari}{2001}]{Saari:2001:chaotic}
	Saari, D.~G. (2001).
	\newblock {\em Chaotic Elections: A Mathematician Looks at Voting}.
	\newblock Providence, RI: American Mathematical Society.
	
	\bibitem[\protect\citeauthoryear{Taylor and Zwicker}{Taylor and
		Zwicker}{1999}]{Taylor/Zwicker:1999}
	Taylor, A.~D. and W.~S. Zwicker (1999).
	\newblock {\em Simple Games}.
	\newblock Princeton, NJ: Princeton University Press.
	
\bibitem[\protect\citeauthoryear{von~Neumann and Morgenstern}{von~Neumann and
  Morgenstern}{1953}]{vonNeumann/Morgenstern:1953}
Von~Neumann, J. and O.~Morgenstern (1953).
\newblock {\em Theory of Games and Economic Behavior\/} (3rd ed.).
\newblock Princeton, NJ: Princeton University Press.
\newblock 
	
	\bibitem[\protect\citeauthoryear{Young}{Young}{1975}]{Young:1975}
	Young, H.~P. (1975).
	\newblock Social choice scoring functions.
	\newblock {\em SIAM Journal on Applied Mathematics\/}~{\em 28\/}(4), 824--838.
	
\end{thebibliography}

%\bibitem[\protect\citeauthoryear{Borda}{Borda}{1784}]{Borda:1784}
%Borda, J.~C. (1784).
%\newblock M\'{e}moire sur les \'{e}lections au scrutin. {{H}istoire et
%  M\'{e}moires de l'Acad\'{e}mie Royale des Sciences, Ann\'{e}e 1781}.
% \newblock Translated by I.~McLean and A.~B. Urken (Eds.), {\em Classics of Social Choice},
%  1995, pp.\  83--89. Ann Arbor, MI: University of Michigan Press.
%

%\bibitem[\protect\citeauthoryear{von~Neumann and Morgenstern}{von~Neumann and
%  Morgenstern}{1953}]{vonNeumann/Morgenstern:1953}
%Von~Neumann, J. and O.~Morgenstern (1953).
%\newblock {\em Theory of Games and Economic Behavior\/} (3rd ed.).
%\newblock Princeton, NJ: Princeton University Press.
%\newblock 

\setstretch{1.24} 

\newpage \small 
\appendix
%\begin{center}
%{\Large\bf -- For online publication --}\vspace{-1cm}
%\end{center}
\section*{Appendix A: Proofs}
    \phantomsection \label{appendix:proofs}
    \addcontentsline{toc}{section}{Appendix A: Proofs}

%\subsection*{Proof of Lemma 1}
%Consider $\mathbf{w}\neq \mathbf{0}\in\mathbb{N}_{0}^n$ and the unanimous profile $\mathbf{P}=(P,\ldots, P)\in\mathcal{P}(A)^n$ with  $a_2 P a_3 P \dots P a_m P a_1$. Then $r|\mathbf{0}(\mathbf{P})=a_1$ but $r|\mathbf{w}(\mathbf{P})=a_2$ for any $r\in \{r^A, r^B, r^C, r^P\}$.

%\subsection*{Proof of Proposition 1} %m=2, all rules coincide
%Let $A=\{a_1,a_2\}$ and fix $\mathbf{w}\neq \mathbf{0}\in \mathbb{N}_0^n$. Then for any $r\in \{r^A, r^B, r^C, r^P\}$
%$$
%r|\mathbf{w}(\mathbf{P})=\begin{cases}
%a_2 & \text{if }\ \sum\limits_{i\colon a_2 P_i a_1}w_i > \sum\limits_{j\colon a_1 P_j a_2}w_j, \\ a_1 & \text{otherwise}.
%\end{cases}
%$$
%%for any $r\in \{r^A, r^B, r^C, r^P\}$. So antiplurality, Borda, Copeland and plurality rule are equivalent and hence have the same equivalence classes.
%
%\subsection*{Proof of Proposition 2} %m=2 = 50% weighted voting
%Define $w(T):=\sum_{i\in T} w_i$ for $T\subseteq N$. If $w(S)\ge w(N\smallsetminus S)$  then $r^P|\mathbf{w}(\mathbf{P})=a_1$ and $v(S)=1$.  
%If $w(S)<w(N\smallsetminus S)$ then $r^P|\mathbf{w}(\mathbf{P})=a_2$ and $v(S)=0$. Proposition~\ref{prop:all_the_same} extends this %observation 
%to $r\in \{r^A, r^B, r^C\}$ . 

\subsection*{Proof of Proposition 3} %antiplurality partitions
The claim is obvious for $n=1$. %, as each non-degenerate weight then is equivalent to $w_1=1$. 
So consider $m\ge n+1$ for $n\ge 2$. 
Let $A=\{a_1, \ldots, a_m\}$ and $\mathbf{P^i}\in\mathcal{P}(A)^n$ be any preference profile where the first $i$ players rank alternative~$a_1$ last and the remaining $n-i$ players rank alternative~$a_2$ last. 
Consider any $\mathbf{\bar w_k}$ and $\mathbf{\bar w_l}$ with $k<l$. Then $r^A|\mathbf{\bar w_k}(\mathbf{P^k})=a_2 \neq r^A|\mathbf{\bar w_l}(\mathbf{P^k})=a_3$. 
So $\mathcal{E}^{r^A}_{\mathbf{\bar w_1},m}, \mathcal{E}^{r^A}_{\mathbf{\bar w_2},m}, 
\ldots, \mathcal{E}^{r^A}_{\mathbf{\bar w_n},m}$ all differ.

Now assume %some 
$\mathbf{w}\in \mathbb{N}_0^n\smallsetminus \{\mathbf{0}\}$ with $ w_1\ge  w_2 \ge \ldots \ge w_n$ satisfies $(r^A,\mathbf{w})\not \sim_m (r^A,\mathbf{\bar w_k})$ for all $k\in \{1, \ldots, n\}$. Let $l$ denote the index such that $w_l>0$ and $w_{l+1}=0$. Then 
both $r^A|\mathbf{w}(\mathbf{P})$ and $r^A|\mathbf{\bar w_l}(\mathbf{P})$ equal the lexicographically minimal element in 
$
Z^l(\mathbf{P}):=\big\{a\in A\ | \ \forall  i\in \{1, \ldots, l\}\colon \exists a'\in A\colon a P_i a' \ \big\}
$, which
collects all alternatives not ranked last by any of the players who have positive weight. These coincide for $\mathbf{w}$ and $\mathbf{\bar w_l}$; and $Z^l(\mathbf{P})$ is non-empty because $m\ge n+1$. This holds for arbitrary $\mathbf{P}\in \mathcal{P}(A)^n$. Hence $r^A|\mathbf{w}\equiv r^A|\mathbf{\bar w_l}$, contradicting the assumption that $(r^A,\mathbf{w})\not \sim_m (r^A,\mathbf{\bar w_k})$ for all $k\in \{1, \ldots, n\}$. Consequently, $\mathcal{E}^{r^A}_{\mathbf{\bar w_1},m}, \mathcal{E}^{r^A}_{\mathbf{\bar w_2},m}, 
\ldots, \mathcal{E}^{r^A}_{\mathbf{\bar w_n},m}$ are all antiplurality classes that exist for $m\ge n+1$ (plus the degenerate $\mathcal{E}_{\mathbf{0},m}$).

\subsection*{Proof of Proposition 4} %Borda partitions
Let $k>j$ for otherwise arbitrary $j, k\in \{1, \ldots, m\}$ 
and consider 
%. Consider $A=\{a_1, \ldots, a_m\}$ and 
any profile $\mathbf{P}\in\mathcal{P}(A)^n$ such that 
%player~1 prefers $a_2$ most and ranks all remaining alternatives lexicographically
%while player~2 ranks $a_2$ in $k$-th position and otherwise agrees with player~1, i.e., suppose 
$a_2\,P_1\,a_1\,P_1\,a_3\,P_1\,a_4\,%P_1\,
\ldots\, %P_1\, 
a_m$ and 
$a_1\,P_2\,a_3\,P_2\,a_4 %\,P_2
\,\ldots$ $%\, P_2\, 
a_k\,P_2\, a_2 \,P_2\, a_{k+1}\,P_2\, a_{k+2}\, %P_2\, 
\ldots \, %P_2\, 
a_{m}.$
The Borda score $j\cdot (m-2) +(m-1)$ of $a_1$ under $\mathbf{\tilde w_j}$ is at least as big as the corresponding score $j\cdot (m-1) + (m-k)$ of $a_2$. Since scores of $a_3, \ldots, a_m$ are all strictly smaller than that of $a_1$, %bounded above by $j\cdot (m-3) +1\cdot (m-2)< j\cdot (m-2) +1\cdot (m-1)$. 
we have $r^B|\mathbf{\tilde w_j}(P)=a_1$. 
With $\mathbf{\tilde w_k}$, by contrast, $a_1$'s weighted score $k\cdot (m-2) + (m-1)$ is strictly smaller than $a_2$'s corresponding score $k\cdot (m-1) + (m-k)$. Scores of $a_3, \ldots, a_m$ remain smaller than $a_1$'s. 
So $r^B|\mathbf{\tilde w_k}(P)=a_2$. Hence $(r^B,\mathbf{\tilde w_j})\not\sim_m (r^B,\mathbf{\tilde w_k})$. 

\subsection*{Proof of Proposition 5} %Copeland partitions 
For a given set %of alternatives 
$A=\{a_1,\dots,a_m\}$ and any subset $A'\subseteq A$ that preserves the order of the alternatives, %we 
denote the \emph{projection} of preference profile $\mathbf{P}\in\mathcal{P}(A)^n$ to $A'$ by $\mathbf{P}\proj{A'}$ 
with $a_k\, P_i\proj{A'} a_l  :\Leftrightarrow [a_k P_i a_l$ and $a_k, a_l\in A']$. 
%For instance, for  $\mathbf{P}=(a_1a_2a_3,a_3a_1a_2,a_2a_3a_1)$ and $A'=\{a_1,a_3\}$ we have $\mathbf{P}\proj{A'}=(a_1a_3,a_3a_1,a_3a_1)$. 
Conversely, if $A'\supseteq A$ is a superset of $A$ with $A'\smallsetminus A=\{a_{m+1}, \ldots, a_{m'}\}$ %we 
define the 
\emph{lifting}  $\mathbf{P}\lift{A'}$ of $\mathbf{P}\in\mathcal{P}(A)^n$ to $A'$ by appending alternatives 
$a_{m+1}, \ldots, a_{m'}$ to each ordering $P_i$ below the lowest-ranked alternative from $A$. 
%That is, for $\mathbf{P}=(a_1a_2a_3,a_3a_1a_2,a_2a_3a_1)$ and $A'=\{a_1,a_2,a_3,a_4\}$ we have $\mathbf{P}\lift{A'}=(a_1a_2a_3 a_4, a_3a_1a_2 a_4, a_2a_3a_1 a_4)$.
%We let $\rho$ or $r$ refer to whole families of mappings and, for instance, write $\rho(\mathbf{P})=\rho(\mathbf{P}\proj{A'})$ if the same alternative $a^*\in A'\subset A$ happens to win for both $A$ and the smaller set $A'$.

Now consider %$A=\{a_1, \ldots, a_m\}$ for 
$m>2$ and %any 
$\mathbf{w},\mathbf{w'}\in\mathbb{N}_{0}^n$ such that $(r^C,\mathbf{w})\not\sim_m (r^C,\mathbf{w'})$, i.e.,
$r^C|\mathbf{w}(\mathbf{P})\neq r^C|\mathbf{w'}(\mathbf{P})$ for some $\mathbf{P}\in\mathcal{P}(A)^n$. 
The $\mathbf{w}$ and $\mathbf{w'}$-weighted versions of the majority relation must differ at $\mathbf{P}$. 
%: if all pairwise comparisons produced the same winners for weights $\mathbf{w}$ and $\mathbf{w'}$, identical Copeland winners would follow.
%i.e., a 
W.l.o.g.\ let weak victory of %some 
$a_k$ over %some 
$a_l$ for $\mathbf{w}$ turn into a strict victory of $a_l$ over $a_k$ for $\mathbf{w'}$:
\setlength{\abovedisplayskip}{8pt}
\setlength{\belowdisplayskip}{8pt}
\be\label{eq:aj_Copelandbeats_ak_and_reverse}
\sum\limits_{i\colon a_k P_i a_l}w_i \ge \sum\limits_{j\colon a_l P_j a_k}w_j
\quad\text{ and }\quad
\sum\limits_{i\colon a_k P_i a_l}w_i' < \sum\limits_{j\colon a_l P_j a_k}w_j'.
\ee
Then take $A'=\{a_k,a_l\}\subset A$ where $|A'|=2$ and projection $\mathbf{P}\proj{A'}$. 
(\ref{eq:aj_Copelandbeats_ak_and_reverse}) implies
\be \label{eq:aj_Copelandbeats_ak_and_reverse2}
\sum\limits_{i\colon a_k \, P_i\ \proj{A'}\, a_l}w_i \ge \sum\limits_{j\colon a_l\, P_j\ \proj{A'}\, a_k}w_j
\quad\text{ and }\quad
\sum\limits_{i\colon a_k \, P_i\ \proj{A'}\,  a_l}w_i' < \sum\limits_{j\colon a_l \, P_j\ \proj{A'}\,  a_k}w_j'.
\ee
If both inequalities are strict or $k<l$ then $r^C|\mathbf{w}(\mathbf{P}\proj{A'})=a_k\neq r^C|\mathbf{w'}(\mathbf{P}\proj{A'})=a_l$ and hence $(r^C,\mathbf{w})\not\sim_2 (r^C,\mathbf{w'})$. 
If not, $a_l$ wins also for $\mathbf{w}$ by lexicographic tie breaking but we can consider profile $\mathbf{P'}\in \mathcal{P}(A')^n$ with $a_l P_i' a_k \Leftrightarrow a_k P_i\proj{A'}\!a_l$ for all $i\in N$. Then
$r^C|\mathbf{w}(\mathbf{P'})=a_l\neq r^C|\mathbf{w'}(\mathbf{P'})=a_k$ and $(r^C,\mathbf{w})\not\sim_2 (r^C,\mathbf{w'})$.

Conversely take $A=\{a_1, a_2\}$ and $\mathbf{w},\mathbf{w'}\in\mathbb{N}_{0}^n$ such that $(r^C,\mathbf{w})\not\sim_2 (r^C,\mathbf{w'})$ and $r^C|\mathbf{w}(\mathbf{P})=a_1\neq r^C|\mathbf{w'}(\mathbf{P})=a_2$ for some $\mathbf{P}\in\mathcal{P}(A)^n$. Then
\be\label{eq:aj_Copelandbeats_ak_and_reverse3}
\sum\limits_{i\colon a_1 P_i a_2}w_i \ge \sum\limits_{j\colon a_2 P_j a_1}w_j
\quad\text{ and }\quad
\sum\limits_{i\colon a_1 P_i a_2}w_i' < \sum\limits_{j\colon a_2 P_j a_1}w_j'.
\ee
Consider $A'=\{a_1,a_2, \ldots, a_m\}\supset A$ where $|A'|=m$ and lifting $\mathbf{P}\lift{A'}$. 
(\ref{eq:aj_Copelandbeats_ak_and_reverse3}) implies
\be 
\sum\limits_{i\colon a_1 \, P_i\ \lift{A'} a_2}w_i \ge \sum\limits_{j\colon a_2\, P_j\ \lift{A'} a_1}w_j
\quad\text{ and }\quad
\sum\limits_{i\colon a_1 \, P_i\ \lift{A'}  a_2}w_i' < \sum\limits_{j\colon a_2 \, P_j\ \lift{A'}  a_1}w_j'
\ee
and alternatives $a_3, \ldots, a_m$ lose all weighted majority comparisons against $a_1$ and $a_2$ by construction of $\mathbf{P}\lift{A'}$. So $r^C|\mathbf{w}(\mathbf{P}\lift{A'})=a_1\neq r^C|\mathbf{w'}(\mathbf{P}\lift{A'})=a_2$. Hence $(r^C,\mathbf{w})\not\sim_m (r^C,\mathbf{w'})$. 
In summary, $(r^C,\mathbf{w})\not\sim_2 (r^C,\mathbf{w'})\Leftrightarrow (r^C,\mathbf{w})\not\sim_m (r^C,\mathbf{w'})$ and, a~fortiori, $(r^C,\mathbf{w})\sim_2 (r^C,\mathbf{w'})\Leftrightarrow (r^C,\mathbf{w})\sim_m (r^C,\mathbf{w'})$.

\subsection*{Proof of Proposition 6} %plurality partitions
Let $m>n$. Consider $A=\{a_1,\dots,a_m\}$ and any  $\mathbf{w},\mathbf{w'}\in\mathbb{N}_{0}^n$ 
such that $(r^P,\mathbf{w})\not\sim_m (r^P,\mathbf{w'})$.
So there exists $\mathbf{P}\in\mathcal{P}(A)^n$ with $r^P|\mathbf{w}(\mathbf{P})=a_k\neq r^P|\mathbf{w'}(\mathbf{P})=a_l$. 
For this $\mathbf{P}$ let
$
\hat A:=\big\{a\ |\  \exists i\in N\colon \forall a'\neq a\colon a P_i a' \big\}
$
denote the set of all alternatives that are top-ranked by some voter. (Obviously, $a_k, a_l\in \hat A$.) Now define $A'\subset A$ as the union of $\hat A$ and some arbitrary elements of $A\smallsetminus \hat A$ such that $|A'|=n$. By construction, each $a\in A'$ has the same weighted number of top positions for projection $\mathbf{P}\proj{A'}$ as it had for $\mathbf{P}$. So $r^P|\mathbf{w}(\mathbf{P}\proj{A'})=a_k\neq r^P|\mathbf{w'}(\mathbf{P}\proj{A'}\nolinebreak)=a_l$. Hence $(r^P,\mathbf{w})\not\sim_n (r^P,\mathbf{w'})$.

Analogously, consider $A=\{a_1,\dots,a_n\}$ and  $\mathbf{w},\mathbf{w'}\in\mathbb{N}_{0}^n$ 
such that $(r^P,\mathbf{w})\not\sim_n (r^P,\mathbf{w'})$. A profile $\mathbf{P}\in\mathcal{P}(A)^n$ with $r^P|\mathbf{w}(\mathbf{P})=a_k\neq r^P|\mathbf{w'}(\mathbf{P})=a_l$ can then be lifted to $A'=A\cup \{a_{n+1},\dots,a_m\}$. By construction, $r^P|\mathbf{w}(\mathbf{P}\lift{A'})=a_k\neq r^P|\mathbf{w'}(\mathbf{P}\lift{A'})=a_l$. Hence $(r^P,\mathbf{w})\not\sim_m (r^P,\mathbf{w'})$. Overall, we can conclude
$(r^P,\mathbf{w})\sim_m (r^P,\mathbf{w'})\Leftrightarrow (r^P,\mathbf{w})\sim_n (r^P,\mathbf{w'})$.

\newpage 

%\begin{center}
%{\Large\bf -- Potentially for ``online only'' publication --}\vspace{-0.8cm}
%\end{center}
\section*{Appendix~B: Minimal representations of committees}
    \phantomsection \label{appendix:representations}
    \addcontentsline{toc}{section}{Appendix~B: Minimal representations}

\vspace{-0.2cm}

\renewcommand{\arraystretch}{1.1}

    \phantomsection 
    \addcontentsline{toc}{subsection}{List of antiplurality games}

\renewcommand{\thetable}{B-\arabic{table}}
\setcounter{table}{0}

\begin{table}[h!] \small
\begin{center}
\begin{tabular}{||c|r c|r c|r c|r c||}
\hline\hline
$n,m$ & \multicolumn{8}{Sc||}{Minimal $\mathbf{\bar w}$ for all antiplurality classes $\mathcal{E}_{\mathbf{\bar w}, m}^{r^A}$ 
}  \\ 
\hline \hline
$3,3$ & 1. & (1,0,0) & 3. & (1,1,1) & 5. & (2,2,1) & &\\
      & 2. & (1,1,0) & 4. & (2,1,1) &     &        & &    \\
\hline\hline
$3, m\ge 4$ & 1. & (1,0,0) & 2. & (1,1,0) & 3. & (1,1,1) & &\\
\hline \hline
$4,3$ & 1. & (1,0,0,0) & 6. & (2,1,1,1) & 11. & (3,2,2,1) & 16. & (4,3,2,2) \\
      & 2. & (1,1,0,0) & 7. & (2,2,1,0) & 12. & (3,3,1,1) & 17. & (4,4,2,1) \\
      & 3. & (1,1,1,0) & 8. & (2,2,1,1) & 13. & (3,3,2,1) & 18. & (4,4,3,2) \\
      & 4. & (1,1,1,1) & 9. & (2,2,2,1) & 14. & (3,3,2,2) & 19. & (5,4,3,2) \\
      & 5. & (2,1,1,0) & 10. & (3,2,1,1) & 15. & (4,3,2,1) &  &  \\
\hline\hline
$4,4$ & 1. & (1,0,0,0) & 3. & (1,1,1,0) & 5. & (2,1,1,1) & 7. & (2,2,2,1) \\
      & 2. & (1,1,0,0) & 4. & (1,1,1,1) & 6. & (2,2,1,1) &  &  \\
\hline
\hline
$4,m\ge 5$ & 1. & (1,0,0,0) & 2. & (1,1,0,0) & 3. & (1,1,1,0) &  4. & (1,1,1,1) \\
\hline
\hline 
\end{tabular}
\caption{Minimal representations of different antiplurality committees\label{table:wcgs_a_3_3}}
\end{center}

\medskip

\subsection*{\ }

    \phantomsection 
    \addcontentsline{toc}{subsection}{List of Borda games}

\begin{center}
\begin{tabular}{||c|rc|rc|rc|rc||}
\hline\hline
$n,m$ & \multicolumn{8}{Sc||}{Minimal $\mathbf{\bar w}$ for all Borda classes $\mathcal{E}_{\mathbf{\bar w}, 3}^{r^B}$} \\
\hline \hline
$3,3$ & 1. & (1,0,0) & 14. & (3,3,2) & 27. & (5,4,3) & 40. & (8,6,3) \\
      & 2. & (1,1,0) & 15. & (4,3,1) & 28. & (7,4,1) & 41. & (9,6,2) \\
      & 3. & (1,1,1) & 16. & (5,2,1) & 29. & (6,5,2) & 42. & (8,7,3) \\
      & 4. & (2,1,0) & 17. & (4,3,2) & 30. & (7,5,1) & 43. & (8,6,5) \\
      & 5. & (2,1,1) & 18. & (5,2,2) & 31. & (6,5,3) & 44. & (10,7,2) \\
      & 6. & (2,2,1) & 19. & (5,3,1) & 32. & (7,5,2) & 45. & (11,7,2) \\
      & 7. & (3,1,1) & 20. & (4,3,3) & 33. & (8,5,1) & 46. & (9,7,5) \\
      & 8. & (3,2,0) & 21. & (5,4,1) & 34. & (6,5,4) & 47. & (10,8,3) \\
      & 9. & (3,2,1) & 22. & (6,3,1) & 35. & (7,5,3) & 48. & (11,8,2) \\
      & 10. & (4,1,1) & 23. & (5,3,3) & 36. & (7,6,2) & 49. & (11,9,3) \\
      & 11. & (3,2,2) & 24. & (5,4,2) & 37. & (8,5,2) & 50. & (13,8,2) \\
      & 12. & (3,3,1) & 25. & (6,4,1) & 38. & (7,5,4) & 51. & (12,9,7) \\
      & 13. & (4,2,1) & 26. & (7,2,2) & 39. & (7,6,4) &  &  \\
\hline\hline
\end{tabular}
\caption{Minimal representations of different Borda committees}
\label{table:wcgs_b_3_3}
\end{center}
\end{table}

\newpage

	\phantomsection 
    \addcontentsline{toc}{subsection}{List of Copeland games}

\begin{table}[h!] \small
\begin{center}
\begin{tabular}{||c|rc|rc|rc|rc||}
\hline\hline
$n$ & \multicolumn{8}{Sc||}{Minimal $\mathbf{\bar w}$ for all Copeland classes $\mathcal{E}_{\mathbf{\bar w}, m}^{r^C}$} \\
&  \multicolumn{8}{Sc||}{and for all classes $\mathcal{E}_{\mathbf{\bar w}, 2}^{r}$ when $r\in \big\{r^A, r^B, r^P\big\}$} \\
& \multicolumn{8}{Sc||}{and for all weighted voting games $[q;\mathbf{w}]$ with $q=0.5\sum w_i$}\\
\hline \hline
3 & 1. & (1,0,0) & 2. & (1,1,0) & 3. & (1,1,1) & 4. & (2,1,1) \\
\hline \hline
4 & 1. & (1,0,0,0) & 4. & (1,1,1,1) & 7. & (2,2,1,1)& &   \\
  & 2. & (1,1,0,0) & 5. & (2,1,1,0) & 8. & (3,1,1,1)& & \\
  & 3. & (1,1,1,0) & 6. & (2,1,1,1) & 9. & (3,2,2,1)& & \\
\hline \hline
5 & 1. & (1,0,0,0,0) & 8. & (2,1,1,1,1) & 15. & (3,2,2,1,0) & 22. & (4,3,2,2,1) \\
  & 2. & (1,1,0,0,0) & 9. & (2,2,1,1,0) & 16. & (4,1,1,1,1) & 23. & (4,3,3,1,1) \\
  & 3. & (1,1,1,0,0) & 10. & (3,1,1,1,0) & 17. & (3,2,2,1,1) & 24. & (5,2,2,2,1) \\
  & 4. & (1,1,1,1,0) & 11. & (2,2,1,1,1) & 18. & (3,2,2,2,1) & 25. & (4,3,3,2,2) \\
  & 5. & (2,1,1,0,0) & 12. & (3,1,1,1,1) & 19. & (3,3,2,1,1) & 26. & (5,3,3,2,1) \\
  & 6. & (1,1,1,1,1) & 13. & (2,2,2,1,1) & 20. & (4,2,2,1,1) & 27. & (5,4,3,2,2) \\
  & 7. & (2,1,1,1,0) & 14. & (3,2,1,1,1) & 21. & (3,3,2,2,2) &  & \\
\hline \hline
6 & \multicolumn{8}{Sc||}{see next page \ldots} \\
\hline %\hline
\end{tabular}
\end{center}
\end{table}

\begin{table}[h!] \small
\begin{center}
\begin{tabular}{||c|rc|rc|rc|rc||}
\hline 
5 & \multicolumn{8}{Sc||}{\ldots see previous page} \\
\hline \hline
6 & 1. & (1,0,0,0,0,0)   & 36. & (3,2,2,2,2,1) & 71. &	(5,4,3,2,1,1) & 106. & (5,5,4,3,3,2) \\
  & 2. & (1,1,0,0,0,0)   & 37. & (3,3,2,2,1,1) & 72. &	(5,4,3,2,2,0) & 107. & (6,4,4,3,3,2) \\
  & 3. & (1,1,1,0,0,0)   & 38. & (3,3,2,2,2,0) & 73. &	(5,4,4,1,1,1) & 108. & (6,5,4,3,2,2) \\
  & 4. & (1,1,1,1,0,0)   & 39. & (3,3,3,1,1,1) & 74. &	(6,3,2,2,2,1) & 109. & (6,5,4,3,3,1) \\
  & 5. & (2,1,1,0,0,0)   & 40. & (4,2,2,2,1,1) & 75. &	(6,3,3,2,1,1) & 110. & (6,5,5,2,2,2) \\
  & 6. & (1,1,1,1,1,0)   & 41. & (4,3,2,1,1,1) & 76. &	(7,2,2,2,2,1) & 111. & (7,4,4,3,2,2) \\
  & 7. & (2,1,1,1,0,0)   & 42. & (4,3,2,2,1,0) & 77. &	(5,4,3,2,2,1) & 112. & (7,5,3,3,2,2) \\
  & 8. & (1,1,1,1,1,1)   & 43. & (4,3,3,1,1,0) & 78. &	(4,4,3,3,2,2) & 113. & (7,5,4,3,2,1) \\
  & 9. & (2,1,1,1,1,0)   & 44. & (5,2,2,1,1,1) & 79. &	(4,4,3,3,3,1) & 114. & (7,5,5,2,2,1) \\
  & 10. & (2,2,1,1,0,0) & 45. &	(5,2,2,2,1,0) & 80. & (5,3,3,3,2,2) & 115. & (8,4,3,3,2,2) \\
  & 11. & (3,1,1,1,0,0) & 46. &	(3,3,2,2,2,1) & 81. & (5,4,3,2,2,2) & 116. & (6,5,4,4,3,2) \\
  & 12. & (2,1,1,1,1,1) & 47. &	(4,3,2,2,1,1) & 82. & (5,4,3,3,2,1) & 117. & (6,5,5,3,3,2) \\
  & 13. & (2,2,1,1,1,0) & 48. &	(4,3,3,1,1,1) & 83. & (5,4,4,2,2,1) & 118. & (7,5,4,3,3,2) \\
  & 14. & (3,1,1,1,1,0) & 49. &	(5,2,2,2,1,1) & 84. & (5,5,3,2,2,1) & 119. & (7,5,4,4,2,2) \\
  & 15. & (2,2,1,1,1,1) & 50. &	(3,3,2,2,2,2) & 85. & (6,3,3,2,2,2) & 120. & (7,5,5,3,3,1) \\
  & 16. & (2,2,2,1,1,0) & 51. &	(3,3,3,2,2,1) & 86. & (6,4,3,2,2,1) & 121. & (7,6,4,3,2,2) \\
  & 17. & (3,1,1,1,1,1) & 52. &	(4,3,2,2,2,1) & 87. & (6,4,3,3,1,1) & 122. & (7,6,4,3,3,1) \\
  & 18. & (3,2,1,1,1,0) & 53. &	(4,3,3,2,1,1) & 88. & (6,4,4,2,1,1) & 123. & (7,6,5,2,2,2) \\
  & 19. & (3,2,2,1,0,0) & 54. &	(4,3,3,2,2,0) & 89. & (7,3,3,2,2,1) & 124. & (8,5,4,3,2,2) \\
  & 20. & (4,1,1,1,1,0) & 55. &	(4,4,2,2,1,1) & 90. & (7,3,3,3,1,1) & 125. & (8,5,5,3,2,1) \\
  & 21. & (2,2,2,1,1,1) & 56.&	(4,4,3,1,1,1) & 91. & (5,4,3,3,3,2) & 126. & (9,4,4,3,2,2) \\
  & 22. & (3,2,1,1,1,1) & 57. &	(5,2,2,2,2,1) & 92. & (5,4,4,3,2,2) & 127. & (7,5,5,4,3,2) \\
  & 23. & (3,2,2,1,1,0) & 58. &	(5,3,2,2,1,1) & 93. & (5,4,4,3,3,1) & 128. & (7,6,5,3,3,2) \\
  & 24. & (4,1,1,1,1,1) & 59. &	(5,3,3,1,1,1) & 94. & (5,5,3,3,3,1) & 129. & (8,5,5,4,2,2) \\
  & 25. & (2,2,2,2,1,1) & 60. &	(5,3,3,2,1,0) & 95. & (5,5,4,2,2,2) & 130. & (8,6,4,3,3,2) \\
  & 26. & (3,2,2,1,1,1) & 61. &	(6,2,2,2,1,1) & 96. & (6,4,3,3,2,2) & 131. & (8,6,5,3,3,1) \\
  & 27. & (3,2,2,2,1,0) & 62. &	(4,3,3,2,2,1) & 97. & (6,4,4,3,2,1) & 132. & (9,5,5,3,2,2) \\
  & 28. & (3,3,1,1,1,1) & 63. &	(5,3,3,2,1,1) & 98. & (6,5,3,2,2,2) & 133. & (7,6,5,4,4,2) \\
  & 29. & (3,3,2,1,1,0) & 64. &	(4,3,3,2,2,2) & 99. & (6,5,3,3,2,1) & 134. & (8,6,5,4,3,2) \\
  & 30. & (4,2,1,1,1,1) & 65. &	(4,3,3,3,2,1) & 100. & (6,5,4,2,2,1) & 135. & (8,7,5,3,3,2) \\
  & 31. & (4,2,2,1,1,0) & 66. &	(4,4,3,2,2,1) & 101. & (7,3,3,3,2,2) & 136. & (9,6,5,4,2,2) \\
  & 32. & (5,1,1,1,1,1) & 67. &	(5,3,2,2,2,2) & 102. & (7,4,3,2,2,2) & 137. & (9,7,5,4,3,2) \\
  & 33. & (3,2,2,2,1,1) & 68. &	(5,3,3,2,2,1) & 103. & (7,4,4,2,2,1) & 138. & (9,7,6,4,4,2) \\
  & 34. & (3,3,2,1,1,1) & 69. &	(5,3,3,3,1,1) & 104. & (7,4,4,3,1,1) &        &                       \\
  & 35. & (4,2,2,1,1,1) & 70. & (5,4,2,2,2,1) & 105. & (8,3,3,3,2,1) &        &                        \\
\hline \hline
\end{tabular}
\end{center}
\caption{Minimal representation of different Copeland committees for $m\ge 2$, \\ and of different antiplurality, Borda and plurality committees for $m=2$,\\ and of different weighted voting games with a simple majority}
\label{table:wcgs_3}
\end{table}

\clearpage

\subsection*{\ }

    \phantomsection 
    \addcontentsline{toc}{subsection}{List of plurality games}

\begin{table} \small 
\begin{center}
\begin{tabular}{||c|rc|rc|rc|rc||}
\hline\hline
$n,m$ & \multicolumn{8}{Sc||}{Minimal $\mathbf{\bar w}$ for all plurality classes $\mathcal{E}_{\mathbf{\bar w}, m}^{r^P}$} \\
\hline \hline
$3, m\ge 3$& 1. & (1,0,0) & 3. & (1,1,1) & 5. & (2,2,1) & & \\
           & 2. & (1,1,0) & 4. & (2,1,1) & 6. & (3,2,2) & & \\
\hline \hline
$4,3$ & 1. & (1,0,0,0) & 10. & (2,2,2,1) & 19. & (4,3,2,1) & 28. & (5,4,3,1) \\
      & 2. & (1,1,0,0) & 11. & (3,2,1,1) & 20. & (4,3,2,2) & 29. & (5,4,3,2) \\
      & 3. & (1,1,1,0) & 12. & (3,2,2,0) & 21. & (4,3,3,1) & 30. & (6,4,3,2) \\
      & 4. & (1,1,1,1) & 13. & (3,2,2,1) & 22. & (4,4,2,1) & 31. & (6,5,3,2) \\
      & 5. & (2,1,1,0) & 14. & (3,3,1,1) & 23. & (5,2,2,2) & 32. & (6,5,4,2) \\
      & 6. & (2,1,1,1) & 15. & (3,2,2,2) & 24. & (4,3,3,2) & 33. & (7,4,4,2) \\
      & 7. & (2,2,1,0) & 16. & (3,3,2,1) & 25. & (5,3,3,1) & 34. & (7,6,4,2) \\
      & 8. & (2,2,1,1) & 17. & (4,2,2,1) & 26. & (5,3,3,2) &  &  \\
      & 9. & (3,1,1,1) & 18. & (3,3,2,2) & 27. & (5,4,2,2) &  & \\
\hline\hline
$4,m\ge 4$ & 1. & (1,0,0,0) & 10. & (2,2,2,1) & 19. & (4,3,2,1) & 28. & (5,4,2,2) \\
        & 2. & (1,1,0,0) & 11. & (3,2,1,1) & 20. & (4,3,2,2) & 29. & (5,4,3,1) \\
        & 3. & (1,1,1,0) & 12. & (3,2,2,0) & 21. & (4,3,3,1) & 30. & (5,4,3,2) \\
        & 4. & (1,1,1,1) & 13. & (3,2,2,1) & 22. & (4,4,2,1) & 31. & (5,4,4,2) \\
        & 5. & (2,1,1,0) & 14. & (3,3,1,1) & 23. & (5,2,2,2) & 32. & (6,4,3,2) \\
        & 6. & (2,1,1,1) & 15. & (3,2,2,2) & 24. & (4,3,3,2) & 33. & (6,5,3,2) \\
        & 7. & (2,2,1,0) & 16. & (3,3,2,1) & 25. & (5,3,3,1) & 34. & (6,5,4,2) \\
        & 8. & (2,2,1,1) & 17. & (4,2,2,1) & 26. & (4,4,3,2) & 35. & (7,4,4,2) \\
        & 9. & (3,1,1,1) & 18. & (3,3,2,2) & 27. & (5,3,3,2) & 36. & (7,6,4,2) \\
\hline\hline
\end{tabular}
\end{center}
\caption{Minimal representations of different plurality committees}
\label{table:wcgs_p_3_3}
\end{table}

\end{document}